\def\DateTime{19/January/1999, 15:30 +0900}
\def\Version{Version 1.0}
\documentclass[leqno,12pt]{amsart}
\usepackage{amssymb}
\usepackage{amscd}
\usepackage{latexsym}
\usepackage{verbatim}

\theoremstyle{plain}
\newtheorem{Theorem}{Theorem}[section]
\newtheorem{Proposition}[Theorem]{Proposition}
\newtheorem{Lemma}[Theorem]{Lemma}
\newtheorem{Corollary}[Theorem]{Corollary}
\newtheorem{Claim}{Claim}[Theorem]

\theoremstyle{definition}

\renewcommand{\theTheorem}{\arabic{section}.\arabic{Theorem}}

\numberwithin{equation}{section}

\newcommand{\ZZ}{{\mathbb{Z}}}
\newcommand{\QQ}{{\mathbb{Q}}}
\newcommand{\RR}{{\mathbb{R}}}
\newcommand{\CC}{{\mathbb{C}}}
\newcommand{\PP}{{\mathbb{P}}}

\newcommand{\FF}{{\mathbb{F}}}
\newcommand{\OO}{{\mathcal{O}}}
\newcommand{\rank}{\operatorname{rk}}
\newcommand{\codim}{\operatorname{codim}}

\newcommand{\GL}{\operatorname{GL}}

\newcommand{\adeg}{\widehat{\operatorname{deg}}}

\newcommand{\Coker}{\operatorname{Coker}}
\newcommand{\Image}{\operatorname{Image}}

\newcommand{\Supp}{\operatorname{Supp}}

\newcommand{\Spec}{\operatorname{Spec}}
\newcommand{\Gal}{\operatorname{Gal}}

\newcommand{\Pic}{\operatorname{Pic}}

\newcommand{\cherncl}{{c}}
\newcommand{\acherncl}{\widehat{{c}}}

\newcommand{\zero}{\operatorname{div}}

\newcommand{\Proof}{{\sl Proof.}\quad}
\newcommand{\QED}{{\unskip\nobreak\hfil\penalty50\quad\null\nobreak\hfil
{$\Box$}\parfillskip0pt\finalhyphendemerits0\par\medskip}}

\newcommand{\Perm}{{\mathfrak{S}}}

\newcommand{\TT}{{\mathcal{T}}}
\newcommand{\HH}{{\mathcal{H}}}
\newcommand{\td}{\operatorname{tr.deg}}
\newcommand{\Bs}{\operatorname{Bs}}

\begin{document}

%%%%%%%%%%%
%% Title %%
%%%%%%%%%%%
\title[Height and intersection]{Height and arithmetic intersection \\ 
       for a family of semi-stable curves}
\author{Shu Kawaguchi}
\address{Department of Mathematics, Faculty of Science,
Kyoto University, Kyoto, 606-01, Japan}
\email[Shu Kawaguchi]{kawaguch@kusm.kyoto-u.ac.jp}
\thanks{Partially supported by JSPS Research Fellowships
for Young Scientists.}
\date{\DateTime, (\Version)}
\begin{abstract}
In this paper, we consider an arithmetic Hodge index theorem 
for a family of semi-stable curves, 
generalizing Faltings-Hriljac's arithmetic Hodge index theorem 
for an arithmetic surface.  
\end{abstract}

\maketitle

%%%%%%%%%%%%%%%%%%%%
%%% Introduction %%%
%%%%%%%%%%%%%%%%%%%%

\section*{Introduction}
\renewcommand{\theTheorem}{\Alph{Theorem}}

In papers \cite{Fa} and \cite{Hr}, Faltings and Hriljac independently 
proved the arithmetic Hodge index theorem on an arithmetic surface. 
Moriwaki \cite{MoH} subsequently proved a higher dimensional case of 
Faltings-Hriljac's arithmetic Hodge index theorem. 
In this paper, we consider an arithmetic Hodge index theorem for 
a family of semi-stable curves. Namely, we prove 

\begin{Theorem}[cf. Theorem~\ref{theorem:main}]
Let $K$ be a finitely generated field over $\QQ$, 
$X_K$ a geometrically irreducible regular projective curve over $K$, 
and $L_K$ a line bundle on $X_K$ with $\deg L_K =0$. 
Let $\overline{B} = (B, \overline{H})$ be a polarization of $K$, 
i.e., 
$B$ a normal projective arithmetic variety 
with the function field $K$, 
and 
$\overline{H}$ a nef $C^{\infty}$-hermitian $\QQ$-line bundle on $B$. 
Let $(X \overset{f}{\longrightarrow} B, \overline{L})$ 
be a model of $(X_K, L_K)$ 
\textup{(}please see \S \textup{4} for terminology\textup{)}. 
We make the following assumptions on the model:
\begin{enumerate}
\renewcommand{\labelenumi}{(\alph{enumi})}
\item
$f$ is semi-stable;  
\item
$X_{\CC}$ and $B_{\CC}$ are non-singular and 
$f_{\CC} : X_{\CC} \to B_{\CC}$ is smooth.  
\end{enumerate}

Let $J_K$ be the Jacobian of $X_K$ and 
$\Theta_{\overline{K}}$ a divisor on 
$J_{\overline{K}}$ which is a translation of 
the theta divisor on $\Pic^{g-1}(X_{\overline{K}})$ 
by a theta characteristic.  
Then we have 
\[
\adeg\left( \acherncl_1(\overline{L})^2 \cdot 
            \acherncl_1\left(f^*(\overline{H})\right)^{d}
     \right)
\leq 
-2 \widehat{h}_{\OO_{J_{\overline{K}}}
(\Theta_{\overline{K}})}^{\overline{B}}([L_K]), 
\]
where $[L_K]$ denotes the point of $J_K$ corresponding to $L_K$  
\textup{(}For the definition of a height function 
$\widehat{h}_{\OO_{J_{\overline{K}}}
(\Theta_{\overline{K}})}^{\overline{B}}$, 
please see \S \textup{4}\textup{)}. 

Furthermore, we assume 
that $H$ is ample and $\cherncl_1(\overline{H})$ is positive. 
Then the equality holds 
if and only if $\overline{L}$ satisfies the following properties:
\begin{enumerate}
\renewcommand{\labelenumi}{(\alph{enumi})}
\item 
There is a Zariski open set $B''$ of $B$ 
with $\codim_B(B \setminus B'') \geq 2$ 
such that $\deg(L\vert_C) = 0$ for any fibral curves $C$ lying over $B''$.
\item 
The restriction of the metric of $\overline{L}$ to each fiber is flat. 
\end{enumerate}
\end{Theorem}

We note that when $B$ is the spectrum of the ring of integers, 
the above theorem is nothing but the arithmetic Hodge index theorem 
for a semi-stable arithmetic surface. 

Our proof uses arithmetic Riemann-Roch theorem, 
similar to that of Faltings on an arithmetic surface, 
although we must consider the Quillen metric. 
Now we outline the organization of this paper. 
In \S 1, we recall some properties of relative Picard 
functors. In \S 2, we recall some facts on determinant 
line bundles, especially for semi-stable curves. 
In \S 3, we deal with an arithmetic setting and 
give hermitian metrics to the results of \S 2. 
In \S 4, we quickly review (a part of) the theory of 
height functions over a finitely generated field over 
$\QQ$, due to Moriwaki \cite{MoA}. 
Finally in \S 5, we prove the main theorem. 

I wish to express my sincere gratitude to Professor Moriwaki 
for his incessant warm encouragement. 
Moreover, it is he who suggested that I consider this work. 

%%%%%%%%%%%%%%%%%%%%%%%%%%
%%% The Picard Functor %%%
%%%%%%%%%%%%%%%%%%%%%%%%%%

\medskip
\section{The Picard functor}
\renewcommand{\theTheorem}{\arabic{section}.\arabic{Theorem}}
The purpose of this section is to review some properties of 
the relative Picard functor, which we will use later. 
We refer to \cite[\S\S 8-9]{Neron} for details. 
In this section, we only deal with schemes which are locally noetherian. 

Let $S$ be a locally noetherian base scheme, $f: X \to S$ 
a flat, projective morphism. 
The {\em relative Picard functor} $\Pic_{X/S}$ of $X$ over $S$ 
is the fppf-sheaf associated with the functor 
\[
P_{X/S}: (\text{locally noetherian $S$-schemes}) 
           \to (\text{Sets}), \quad
T \mapsto \Pic(X \times_S T). 
\]

If we assume $f_*(\OO_X)=\OO_S$ holds universally, 
then for all locally noetherian $S$-schemes $g : T \to S$, 
\[
\Pic_{X/S}(T) = \Pic(X \times_S T)/\Pic(T). 
\]

Furthermore, if $X/S$ admits a section $\epsilon : S \to X$, 
then one checks immediately, 
\begin{equation}
\label{eqn:pic}
\Pic_{X/S}(T) 
= \left\{
  \begin{gathered}
       \text{group of isomorphism classes of} \\                   
       \text{invertible sheaves $L$ on $X \times_S T$,} \\
       \text{plus isomorphism 
          $(\epsilon \circ g, 1_T)^*(L) \simeq \OO_T$}
  \end{gathered}
  \right\}.
\end{equation}
Such invertible sheaves are said to be {\em rigidified} along 
the induced section $\epsilon_T = \epsilon \circ g$. 

If $S$ consists of a field, then $\Pic_{X/S}$ is a group scheme. 
Let $\Pic^{0}_{X/S}$ be its identity component. 
For a general locally noetherian scheme $S$, 
we introduce $\Pic^{0}_{X/S}$ as the subfunctor of $\Pic_{X/S}$ 
which consists of all elements whose restrictions to all fibers $X_s$, 
$s$ being a point of $S$, belong to $\Pic^{0}_{X_{s}/k(s)}$. 

If $X$ is a proper curve over a field $k$, 
then $\Pic^0_{X/k}$ consists of all elements of $\Pic_{X/k}$ 
whose partial degree on each irreducible components of 
$X \otimes_k \overline{k}$ is zero, 
where $\overline{k}$ is an algebraic closure of $k$. 

We note that 
if $\Pic_{X/S}$ (resp. $\Pic^{0}_{X/S}$) is representable 
by a locally noetherian scheme, 
then for all locally noetherian $S$-schemes $T$, 
\[
\Pic_{X/S} \times_S T = \Pic_{X \times_S T/ T} \quad 
(\text{resp. } \Pic^0_{X/S} \times_S T = \Pic^0_{X \times_S T/ T}). 
\]

Now we introduce the notion of universal line bundles 
when $\Pic_{X/S}$ (resp. $\Pic^{0}_{X/S}$) is representable 
by a locally noetherian scheme. 
We assume that 
the structural morphism $f: X \to S$ admits a section $\epsilon$ 
and that $f_*(\OO_X)=\OO_S$ holds universally, 
so that $\Pic_{X/S}$ is given by \eqref{eqn:pic} 
for a locally noetherian $S$-scheme. 
If $\Pic_{X/S}$ (resp. $\Pic^{0}_{X/S}$) is representable 
by a locally noetherian scheme, then the identity on 
$\Pic_{X/S}$ (resp. $\Pic^{0}_{X/S}$) gives rise to a line bundle 
$U$ (resp. $U^0$) 
on $X \times_S \Pic_{X/S}$ (resp. $X \times_S \Pic^{0}_{X/S}$) 
which is rigidified along the induced section. 
$U$ (resp. $U^0$) is called the 
{\em universal} line bundle. 
The justification of the notion of ``universal'' is: 

\begin{Proposition}
\label{prop:universal}
Let $f: X \to S$ be a flat morphism of locally noetherian schemes  
and let $\epsilon$ be a section of $f$. 
Assume that $f_*(\OO_X)=\OO_S$ holds universally. 
If $\Pic_{X/S}$ (resp. $\Pic^{0}_{X/S}$) is representable 
by a locally noetherian scheme, 
the universal line bundle $U$ has the following property: 
For every locally noetherian scheme $g : T \to S$, 
and for every line bundle $L'$ on $X' = X \times_S T$ 
which is rigidified along the induced section 
$\epsilon' = \epsilon \circ g$, 
there exists a unique morphism 
$g: T \to \Pic_{X/S}$ such that $L'$ 
is isomorphic to $(1 \times g)^*(U)$. 

If $\Pic^{0}_{X/S}$ is representable 
by a locally noetherian scheme, 
the universal line bundle $U^{0}$ has a similar property 
for a line bundle $L'$ on $X' = X \times_S T$ 
which is rigidified along the induced section 
and $L_{t}' \in \Pic^{0}_{X_t/k(t)}$ for all $t \in T$.  
\end{Proposition}

\Proof 
\cite[8.2. Proposition 4]{Neron}
\QED

Now we restrict ourselves to the case of semi-stable curves. 
We recall that a {\em semi-stable curve of genus $g$} is 
a proper flat morphism $f : X \to S$ whose fiber 
$X_{\overline{s}}$ over every geometric point $\overline{s}$ of $S$ 
is a reduced connected curve 
with at most ordinary double points 
such that $\dim_{k(s)} H^1(X_{\overline{s}},\OO_{X_{\overline{s}}})$ 
equals to $g$. 

\begin{Proposition}
Let $f: X \to S$ be a semi-stable curve of locally noetherian schemes. 
Then $f_*(\OO_X) = \OO_S$ holds universally. 
\end{Proposition}

\Proof
We have only to prove that $f_*(\OO_X) = \OO_S$. 
Let $\pi \circ \widetilde{f}$ be the Stein factorization of $f$, 
where $\widetilde{f} : X \to \widetilde{S}$ is a proper morphism 
with connected fibers and 
$\pi : \widetilde{S} \to S$ is a finite morphism. 
Since every fiber is geometrically reduced and 
geometrically connected, 
there is a section $\eta : S \to \widetilde{S}$ 
such that $\widetilde{f} = \eta \circ f$ 
by rigidity lemma (\cite[Proposition~6.1]{GIT}).  
Since $\OO_{\widetilde{S}} \simeq \widetilde{f}_*(\OO_{X})$ 
factors through
\[ 
\OO_{\widetilde{S}} \to \eta_*(\OO_S) 
\to f_* \eta_* (\OO_S) = \widetilde{f}_*(\OO_X),  
\]
$\OO_{\widetilde{S}} \to \eta_*(\OO_S)$ is injective. 
On the other hand, 
since $\eta$ is a closed immersion, 
$\OO_{\widetilde{S}} \to \eta_*(\OO_S)$ is surjective, 
hence $\OO_{\widetilde{S}} = \eta_*(\OO_S)$. 
Then, $f_*(\OO_X) = \pi_* \widetilde{f}_*(\OO_X)
= \pi_*(\OO_{\widetilde{S}}) 
= \pi_*(\eta_*(\OO_S)) = \OO_S$.  
\QED

We finish this section by quoting a result obtained by Deligne 
concerning the representability of the relative Picard functor. 

\begin{Theorem}
\label{theorem:Deligne}
Let $f: X \to S$ be a semi-stable curve of locally noetherian schemes. 
Then $\Pic_{X/S}$ is a smooth algebraic space over $S$. 
The identity component $\Pic^{0}_{X/S}$ is a semi-abelian scheme. 
\end{Theorem}

\Proof
\cite[9.4. Theorem 1]{Neron} or 
\cite[Proposition 4.3]{De}
\QED

\medskip

%%%%%%%%%%%%%%%%%%%%%%%%%%%%%%%%%%%
%%% The determinant line bundle %%%
%%%%%%%%%%%%%%%%%%%%%%%%%%%%%%%%%%%

\section{Determinant line bundles}
The purpose of this section is to review some properties of 
determinant line bundles. 
Since we are concerned about a family of curves in this paper, 
we only consider determinant line bundles in a restricted context. 
For a general treatment of determinant line bundles, 
we refer to \cite{KM}.

\begin{Theorem}
\label{theorem:det}
Let us consider a morphism $f: X \to S$ of noetherian schemes 
with the following conditions: 
\begin{enumerate}
\renewcommand{\labelenumi}{(\roman{enumi})}
\item
$f$ is proper, $f_*(\OO_X)=\OO_S$, and $\dim f =1$;
\item
There is an effective Cartier divisor $D$ on $X$ such that 
$D$ is $f$-ample and flat over $S$. 
\end{enumerate}
For every $f: X \to S$ satisfying the above conditions, 
for every line bundle $L$ on $X$ 
and isomorphism of sheaves 
$\phi: L \overset{\sim}{\longrightarrow} L'$, 
one can uniquely construct 
a line bundle $\det Rf_*(L)$ on $S$ and an isomorphism 
$\det Rf_*(\phi): 
\det Rf_*(L) \overset{\sim}{\longrightarrow} \det Rf_*(L')$ 
in such a way that $\det Rf_*(L)$ becomes a functor 
with the following properties: 
\begin{enumerate}
\renewcommand{\labelenumi}{(\alph{enumi})}
\item
If $f_*(L)$ and $R^{1}f_*(L)$ are both locally free, then 
\[
\det Rf_*(L) = \det f_*(L) \otimes 
               \left( \det R^{1}f_*(L) \right)^{-1}; 
\]
\item
$\det Rf_*(L)$ is compatible with a base change, i.e., 
if $g : T \to S$ is a morphism of noetherian schemes, then
\[
g^*\left( \det Rf_*(L) \right) \cong 
\det R(f_T)_*(L_T);
\] 
\item
If $S$ is connected and $M$ is a line bundle on $S$, then 
\[
\det Rf_*\left(L \otimes f^*(M)\right)
\cong \det Rf_*(L) \otimes M^{\chi}, 
\]
where $\chi = \chi(C_s,L_s)$ for some $s \in S$; 
\item
If $D$ is an effective Cartier divisor on $X$ 
which is flat over $S$, then
\[
\det Rf_*(L) 
\cong \det Rf_*\left(L(-D)\right) \otimes \det f_*(L \vert_D).
\]
\end{enumerate}
\end{Theorem}

\Proof
\cite{KM} or \cite[VI \S 6]{La}
\QED

Suppose now that 
$f: X \to S$ is a semi-stable curve of noetherian schemes 
and assume that 
$f$ admits a section $\epsilon$.  
Moreover, let $A$ be a rigidified line bundle on $X$ of  degree $g-1$. 
By Theorem~\ref{theorem:Deligne}, 
$\Pic^{0}_{X/S}$ is a semi-abelian scheme and 
there exists a universal line bundle $U^0$ on $X \times_S \Pic^0_{X/S}$. 
Let $P^a$ be the scheme 
which is the translation of $\Pic^{0}_{X/S}$ by $A$, i.e., 
\begin{equation*}
P^a(T) 
= \left\{
  \begin{gathered} 
          \text{rigidified line bundle $L$ on $X_T$} \\ 
          \text{such that $L \otimes A^{-1}$ 
                    belongs to $\Pic^{0}_{X/S}$}
  \end{gathered} 
  \right\}.
\end{equation*}
Moreover, let $U^a$ be the line bundle on $P^a$ 
which is the translation of $U^0$ by $A$. 
If $q^a : X \times_S P^a \to P^a$ is the second projection, 
then $q^a$ satisfies the condition of Theorem~\ref{theorem:det}, 
because $f : X \to S$ satisfies the condition of Theorem~\ref{theorem:det}.  
Thus the determinant line bundle $\det Rq^a_*(U^a)$ on $P^a$ is defined. 
To simplify the notation, let us denote $\det Rq^a_*(U^a)$ by $\TT^{-1}$. 

In the following, 
we will see that $\TT^{-1}$ is related to the theta divisor. 
Here we further assume that $f: X \to S$ is smooth of genus $g \geq 1$. 
First, we define the theta divisor. 

Let $(X/S)^{(g-1)}$ be the symmetric $(g-1)$-fold product, i.e., 
\[
(X/S)^{(g-1)} = 
\overbrace{X \times_S \cdots \times_S X}^{\text{$(g-1)$ times}} 
/ \Perm_{g-1}, 
\]
where the $(g-1)$-th symmetric group $\Perm_{g-1}$ acts on 
$\overbrace{X \times_S \cdots \times_S X}^{\text{$(g-1)$ times}}$ 
naturally. 
Let 
\[
(X/S)^{(g-1)} \to \Pic^{g-1}_{X/S}, \quad D_T \to [D_T]
\]
be a morphism, 
where for any locally noetherian $S$-schemes $T$ 
and for any $T$-valued point $D_T$ of $(X/S)^{(g-1)}$ 
(i.e., for any effective Cartier divisors on $X \times_S T$ of 
degree $(g-1)$), 
we denote by $[D_T]$ the element of $\Pic^{g-1}_{X/S}$ 
corresponding to $D_T$. 
The schematic image of this morphism, 
which turns out an effective relative Cartier divisor 
on $\Pic^{g-1}_{X/S}$, is called the theta divisor for 
$X/S$ and denoted by $\Theta_{X/S}$. 

\begin{Proposition}
\label{prop:theta}
Let $f : X \to S$ be a projective smooth morphism of noetherian schemes 
whose geometric fibers are smooth projective curves of genus $g \geq 1$. 
We assume the existence of a section. 
Let $\Pic^{g-1}_{X/S}$ be a Picard scheme of degree $(g-1)$ 
and $U$ a universal line bundle on $X \times_S \Pic^{g-1}_{X/S}$. 
Then 
\[
\det Rq_*(U) \cong \OO_{\Pic^{g-1}_{X/S}}(-\Theta_{X/S}),
\]
where $\Theta_{X/S}$ is the theta divisor for $X/S$ and 
$q: X \times_S \Pic^{g-1}_{X/S} \to \Pic^{g-1}_{X/S}$ is 
the second projection. 
\end{Proposition}

\Proof
When the base scheme is a point, or an arithmetic surface, 
this is well-known 
(cf. \cite[\S 5]{Fa} or \cite[VI Lemma~2.4]{La}). 
The proof for a general base scheme is similar to 
that for a point, as we will see in the following. 

Let $p : X \times_S \Pic^{g-1}_{X/S} \to \Pic^{g-1}_{X/S}$ 
be the first projection. 
Let $D'$ be an effective relative Cartier divisor of 
sufficiently large degree on $X$ 
(actually $\deg D' \geq g$ is enough) 
and put $D = p^*(D')$. 
Since 
\[
H^0(X_s, U(-D)_t) =0
\]
for all points $t$ of $\Pic^{g-1}_{X/S}$ 
and the point $s$ of $S$ lying below $t$, 
$q_*(U(-D)) =0$ by \cite[Corollorary II.12.9]{Ha}, 
and $R^1q_*(U(-D))$ is locally free. 
Thus, by (a) and (d) of Theorem~\ref{theorem:det}, 
\[
\det Rq_*(U) = \det q_*(U\vert_{D}) 
                 \otimes \left( R^1q_*(U(-D)) \right)^{-1}. 
\]
Since $q_*(U)$ is torsion-free and 
$H^0(X_s,U_t) =0$ for a general point $t$ of $P$, 
it follows that $q_*(U)=0$. 
Also, since $D \to \Pic^{g-1}_{X/S}$ is finite, 
$R^1q_*(U\vert_{D})=0$. 
Thus we get the exact sequence:
\[
0 \to q_*(U\vert_{D}) \to R^1q_*(U(-D)) \to R^1q_*(U) \to 0. 
\]
We denote the homomorphism 
$q_*(U\vert_{D}) \to R^1q_*(U(-D))$ by $\alpha$. 
Since $R^2q_*(U) =0$, we get by \cite[Theorem II.12.11]{Ha}
\[
R^1q_*(U) \otimes k(t) \cong H^1(X_s,U_t)
\]
for all points $t$ of $\Pic^{g-1}_{X/S}$ 
and the point $s$ of $S$ lying below $t$. 
If $R^1q_*(U) \otimes k(t) =0$, 
then $R^1q_*(U)$ is also zero for some neighborhood of $t$, 
and especially $R^1q_*(U)$ is flat for some neighborhood 
of $t$. Thus
\begin{align*}
\text{$\alpha(t)$ is an isomorphism}  
  & \Leftrightarrow R^1q_*(U) \otimes k(t) =0  \\ 
  & \Leftrightarrow H^1(X_s,U_t)=0 \\ 
  & \Leftrightarrow t \not\in \Theta_{X/S}. 
\end{align*}
Therefore if we put 
$E = \{ t \in \Pic^{g-1}_{X/S} \mid (\det\alpha)(t) = 0\}$, 
then $E = a\Theta_{X/S}$ for some positive integer $a$. 
By considering the case that the base scheme is a point, 
we get $a =1$. 
\QED

Now we put everything together and get: 

\begin{Theorem}
\label{theorem:det:for:semistable}
Let $f : X \to S$ be a semi-stable curve of genus $g \geq 1$
of noetherian schemes 
and assume that $f$ admits a section $\epsilon$. 
Let $A$ be a rigidified line bundle of degree $(g-1)$ 
and $(P^a, U^a)$ the translation of $(\Pic^0_{X/S}, U^0)$ by $A$.   
We put $\TT^{-1} = \det Rq^a_*(U^a)$, 
where $q^a : X \times_S P^a \to P^a$ is the second projection. 
Then, 
\begin{enumerate}
\renewcommand{\labelenumi}{(\roman{enumi})}
\item
If $T \to S$ be a morphism of noetherian schemes such that 
$f_T : X \times_S T \to T$ is smooth, then 
\[
\TT_{T}^{-1} = \OO_{P^a_T}(-\Theta_{X_T/T})
\]
where $\Theta_{X_T/T}$ is the theta divisor for $X_T/T$. 
\item
If $L$ is a rigidified line bundle on $X$ 
which belongs to $P^a(S)$, 
then there is a canonical morphism $g^a : S \to P^a$ 
such that the induced morphism 
\[
u_L : \det Rf_*(L) \to (g^a)^*(\TT^{-1})
\] 
is canonically isomorphic.  
\end{enumerate} 
\end{Theorem}

\Proof
Noting that determinant line bundles are compatible with a base change, 
we have already seen (i).
Regarding as (ii), by the universal property of $U^a$, 
there exists a canonical morphism $g^a : S \to P^a$ such that 
\[
L \cong (1 \times g^a)^*(U^a). 
\]
On the other hand, 
since determinant line bundles are compatible with a base change, 
we have canonically 
\[
(g^a)^* \left( \det Rq^a_*(U^a) \right) 
\cong \det Rf_*(\left (1 \times g^a)^*(U^a)\right). 
\] 
Combining above two isomorphisms, 
we get the desired isomorphism.  
\QED

%%%%%%%%%%%%%%%%%%%%%%%%%%%%
%%% Arithmetic Situation %%%
%%%%%%%%%%%%%%%%%%%%%%%%%%%%

\medskip
\section{Arithmetic Setting}
In this section, we consider an arithmetic setting. 
An {\em arithmetic variety}  
is an integral scheme  
which is flat and quasi-projective over $\Spec(\ZZ)$. 

Let $f: X \to B$ be a semi-stable curve of genus $g \geq 1$ 
of arithmetic varieties   
and assume that
$f$ admits a section $\epsilon$. 
We also assume that $f_{\CC} : X_{\CC} \to B_{\CC}$ 
is a smooth morphism. 
Let $A$ be a rigidified line bundle of degree $(g-1)$ and 
$(P^a, U^a)$ the translation of 
$(\Pic^0_{X/S}, U^0)$ by $A$. 
We put $\TT^{-1} = \det Rq^a_*(U^a)$ on $P^a$, 
where $q^a : X \times P^a \to P^a$ is the second projection. 
Then by Theorem~\ref{theorem:det:for:semistable}(ii), 
for a rigidified line bundle $L$ which belongs to $\Pic^0_{X/S}$, 
we have a natural isomorphism 
\[
u_L : \det Rf_*(L \otimes A) \to (g^a)^*(\TT^{-1}), 
\]
where $g^a : S \to P^a$ is an induced morphism by $L \otimes A$. 

In this section we give metrization on the above line bundles, 
and consider the norm of $u_L$. 
Let $\Theta_{X_{\CC}/B_{\CC}}$ 
be the theta divisor for $X_{\CC}/B_{\CC}$, 
which is a relative Cartier divisor on 
$P^a_{\CC} = \Pic^{g-1}_{X_{\CC}/B_{\CC}}$. 
Then by Theorem~\ref{theorem:det:for:semistable}(i), 
$\TT^{-1}_{\CC} 
= \OO_{\Pic^{g-1}_{X_{\CC}/B_{\CC}}}(-\Theta_{X_{\CC}/B_{\CC}})$.

In the following, we introduce a metric on 
$\OO_{\Pic^{g-1}_{X_{\CC}/B_{\CC}}}(-\Theta_{X_{\CC}/B_{\CC}})$. 
Put $J = \Pic^{0}_{X_{\CC}/B_{\CC}}$ and let 
\[
\lambda : \Pic^{g-1}_{X_{\CC}/B_{\CC}} \to J \quad
[D_T] \mapsto (g-1)[\epsilon_T]
\]
be an isomorphism, where for any $B_{\CC}$-scheme $T$, 
$[\epsilon_T]$ is the class of the induced section by $\epsilon$. 
Let $\Theta^0_{X_{\CC}/B_{\CC}}$ be the image of 
$\Theta_{X_{\CC}/B_{\CC}}$ by $\lambda$.   

We need some definitions to proceed. 
The {\em Siegel upper-half space of degree $g$}, 
denoted by $\HH_g$, is defined by 
\[
\HH_g = \{ \Omega = X + \sqrt{-1} Y \in \GL_g(\CC) \mid 
{}^t\Omega = \Omega, Y > 0 \}. 
\]
Moreover, the {\em symplectic group of degree $2g$}, 
denoted by $Sp_g(\ZZ)$, is defined by 
\[
Sp_g(\ZZ) = \{ S \in GL_{2g}(\ZZ) \mid {}^t S J S = J \}, 
\]
where $J = \begin{pmatrix}
             0 & -I \\
            -I & 0 \\
           \end{pmatrix}$.  
An element $S =  \begin{pmatrix}
                   A & B \\
                   C & D \\
                  \end{pmatrix}$ 
of $Sp_g(\ZZ)$ acts on $\HH_g$ by 
\[
S \cdot \Omega = (A \Omega + B) (C \Omega +D)^{-1} 
\]
and $Sp_g(\ZZ) \backslash \HH_g$ becomes a coarse moduli of 
principally polarized abelian varieties. 

For $z = x + \sqrt{-1}y \in \CC^g$ and 
$\Omega = X + \sqrt{-1}Y \in \HH_g$, 
we define 
\begin{gather*}
\theta(z,\Omega) 
= \sum_{m \in \ZZ^g} \exp (\pi \sqrt{-1} \ {}^t m\Omega m 
  + 2\pi\sqrt{-1} \ {}^t m \cdot z), \\
\Vert \theta \Vert (z,\Omega)
= \sqrt[4]{\det Y} \exp(-\pi {}^t y Y y) \vert \theta(z,\Omega) \vert.  
\end{gather*}
Then $\theta$ becomes a holomorphic function on $\CC^g \times \HH_g$. 
Moreover $\Vert \theta \Vert$ becomes a $C^{\infty}$-function 
which is periodic with respect to $\ZZ^g + \Omega \ZZ^g$, 
so that $\Vert \theta \Vert$ is seen as a $C^{\infty}$-function 
on $\CC^g/ \ZZ^g + \Omega \ZZ^g$. 

Going back to our situations, 
for any $b \in B(\CC)$, let us write analytically 
\[
J_b \cong \CC^g/ \ZZ^g + \Omega_b \ZZ^g
\]
where $\Omega_b \in \HH_g$. 
Then there is a unique element $t_b \in \CC^g/ \ZZ^g + \Omega \ZZ^g$ 
such that $\Theta^0_{X_b} = \zero \left( \theta(z +t_b, \Omega_b) \right)$, 
where $\theta(z +t_b, \Omega_b)$ is seen as a function of $z$. 

\begin{Proposition}
\label{prop:metric:of:theta}
Let the notation be as above. 
Let $1$ denote the section of 
$\OO_{J}(\Theta^0_{X_{\CC}/B_{\CC}})$ 
which corresponds to $\Theta^0_{X_{\CC}/B_{\CC}}$. 
For any $p \in J$, 
let $b \in B(\CC)$ be the point lying below $p$ and 
write $J_b \cong \CC^g/ \ZZ^g + \Omega_b \ZZ^g$ and 
$\Theta^0_{X_b} = \zero \left( \theta(z +t_b, \Omega_b) \right)$ 
with $t_b \in \CC^g/ \ZZ^g + \Omega \ZZ^g$. 
Moreover, let $z \in \CC^g/ \ZZ^g + \Omega \ZZ^g$ correspond to $p$. 
Then, if we define 
\[
\Vert 1 \Vert_{\Theta^0_{X_{\CC}/B_{\CC}}} (p)
= 
\Vert \theta \Vert(z + t_b, \Omega_b), 
\]
then $\Vert \cdot \Vert_{\Theta^0_{X_{\CC}/B_{\CC}}}$ 
gives a $C^{\infty}$ metric on 
$\OO_{J}(\Theta^0_{X_{\CC}/B_{\CC}})$ 
\end{Proposition}

\Proof
If the base space $B(\CC)$ is a point, 
the assertion is well-known (cf. \cite[\S 3]{Fa}). 
Thus all we need to prove is that 
$\Vert 1 \Vert_{\Theta^0_{X_{\CC}/B_{\CC}}}$ varies 
smoothly as $b \in B(\CC)$ varies. 
However, since the morphism 
\[
\Phi : B(\CC) \to Sp_g(\ZZ) \backslash \HH_g, \quad
b \mapsto \text{the class of $J_b$}
\]
is holomorphic and 
$t_b$ is given the difference of the section 
$\epsilon_{\CC}$ and a theta characteristic, 
$\Vert 1 \Vert_{\Theta^0_{X_{\CC}/B_{\CC}}}$ varies 
smoothly as $b \in B(\CC)$ varies. 
\QED

Finally, $\OO_{\Pic^{g-1}_{X_{\CC}/B_{\CC}}}(\Theta_{X_{\CC}/B_{\CC}})$ 
is metrized by 
$\left(\OO_{J}(\Theta^0_{X_{\CC}/B_{\CC}}),   
\Vert \cdot \Vert_{\Theta^0_{X_{\CC}/B_{\CC}}} \right)$ 
through $\lambda$. 
We write this metric 
by $\Vert \cdot \Vert_{\Theta_{X_{\CC}/B_{\CC}}}$. 

\medskip
Next we give a $C^{\infty}$ metric on $L_{\CC}$ over $X_{\CC}$. 
Actually, there is a certain class of $C^{\infty}$ metrics on $L_{\CC}$ 
which is suitable for our purpose. 
We introduce this class in the following. 

First we recall admissible metrics of line bundles  
on a compact Riemann surface. 
Let $M$ be a compact Riemann surface of genus $g \geq 1$ and 
$\{ \omega_1, \omega_2, \cdots, \omega_g \}$ 
a basis of $H^0(M, \Omega^1_M)$ with 
\[
\frac{\sqrt{-1}}{2} \int_M \omega_i \wedge \overline{\omega_j} = \delta_{ij}.
\] 
Let us put 
\[
\mu = \frac{\sqrt{-1}}{2g} \sum_{i=1}^{g} \omega_i \wedge \overline{\omega_i}. 
\]
Then $\mu$ is a positive $(1,1)$-form on $M$, 
and is called the canonical volume form on $M$.  
A $C^{\infty}$-metric $h_L$ of a line bundle $L$ on $M$ 
is said to be {\em admissible} if 
\[
\cherncl_1\left((L, h_L)\right) = (\deg L) \mu.
\]   
For every line bundle on $M$, we can endow an admissible metric 
unique up to a constant multiplication. 

Now let us go back to our situation, 
i.e., the case that $f : X_{\CC} \to B_{\CC}$ is a smooth family of 
curves of genus $g \geq 1$.  
A $C^{\infty}$-metric $h_L$ on $L_{\CC}$ over $X_{\CC}$ 
is said to be {\em admissible} 
if for any $b \in B(\CC)$, 
its restriction $(L_b, h_{L,b})$ on $X_b$ is admissible. 
The existence of an admissible metric is given by: 

\begin{Proposition}
\label{prop:admissible:metric}
Let $X$ and $B$ be smooth varieties over $\CC$ 
and $f : X \to B$ a smooth projective morphism with a section 
whose fibers are curves of genus $g \geq 1$. 
Let $L$ be a line bundle on $X$. 
Then there exists a \textup{(}global\textup{)} 
admissible metric on $L$ over $X$. 
\end{Proposition}

\Proof
First we construct a suitable $(1,1)$-form on $X$. 
Let 
\[
j : X \to J = \Pic^0_{X/B} 
\]
is the embedding induced by the section. 
On $J$, we have a $C^{\infty}$-hermitian line bundle 
$\left( \OO_J(\Theta^0_{X/B}), \Vert \cdot \Vert_{\Theta^0_{X/B}}\right)$ 
by Proposition~\ref{prop:metric:of:theta}. 
We consider 
\[
\omega = \frac{1}{g} 
  j^* \left( \cherncl_1 
        \left(  \OO_J(\Theta^0_{X/B}), \Vert \cdot \Vert_{\Theta^0_{X/B}}
        \right)
      \right). 
\]
Then, for any $b \in B$, $\omega_b = \omega\vert_{X_b}$ 
is the canonical volume form on $X_b$ (cf. \cite[Thoerem 1]{Fa}). 

Next we note that the statement of the proposition 
is local with respect to $B$. 
Indeed, let $\{U_i\}_{i=0}^{\infty}$ be an open covering of $B$ 
and assume that for each $i$ there is an admissible metric 
$h_L^{(i)}$ on $L\vert_{f^{-1}(U_i)}$ over $f^{-1}(U_i)$. 
Take a partition of unity $\{\rho_i\}$ subordinate to $\{U_i\}$. 
Then it is easy to see that 
\[
\sum_{i=0}^{\infty} f^*(\rho_i) h_L^{(i)}
\]
is an admissible metric on $X$. 

For $b \in B$, we consider a small ball $U \subset B$ containing $b$. 
We set $X_U = f^{-1}(U)$. 
Then there is a diffeomorphism 
$g : X_U \overset{\approx}{\longrightarrow} X_b \times U$ over $U$ 
(\cite[Theorem~2.4]{Ko}). 

Let $h_L$ be any $C^{\infty}$-hermitian metric on $L$ over $X_U$ 
and set 
$\eta = \cherncl_1(L,h_L)$. 
Then $\eta$ is a $d$-closed $(1,1)$-form. 
We claim that $\left( \deg(L)\omega - \eta \right) \vert_{X_U}$ 
is $d$-exact. 
Indeed, $\left( \deg(L)\omega - \eta \right) \vert_{X_b} = 0$ 
in $H^2(X_b,\CC)$. 
On the other hand, $H^2(X_U,\CC) = H^2(X_b,\CC)$ by Poincar\'{e}'s lemma. 
Thus $\left( \deg(L)\omega - \eta \right) \vert_{X_U}$ 
is a $d$-exact $(1,1)$-form. 
Then by the $d d^c$-lemma, there is a $C^{\infty}$ function $\psi$ 
on $X_U$ with 
$d d^c (\psi) = \left( \deg(L)\omega - \eta \right) \vert_{X_U}$. 
Now if we set 
$h_{L,U}' = \exp(-\psi) h_L\vert_{X_U}$ on $L\vert_{X_U}$ over $X_U$, 
then 
\[
d d^c \left( \cherncl_1(L\vert_{X_U},h_{L,U}') \right) 
= \deg(L) \omega \vert_{X_U},
\]
which is an admissible metric on $L\vert_{X_U}$ over $X_U$.  
\QED

Now we prove the main proposition of this section, 
which will be a key point to prove 
Proposition~\ref{prop:main}.  

\begin{Proposition}
\label{prop:key:lemma}
Let $f: X \to B$ be a semi-stable curve of genus $g \geq 1$ 
of arithmetic varieties  
and assume that
$f$ admits a section $\epsilon$. 
We also assume that $f_{\CC} : X_{\CC} \to B_{\CC}$ 
is a smooth morphism. Let
\[
u_L : \det Rf_*(L \otimes A) \to (g^a)^*(\TT^{-1}), 
\]
be the isomorphism given at the beginning of this section. 
We endow $C^{\infty}$ metrics on $A$ and $\omega_{X/B}$, 
an admissible metric on $L$, 
and then the Quillen metric on $\det Rf_*(L \otimes A)$ 
determined by these metrics. 
Moreover, we endow a metric 
$\Vert \cdot \Vert_{\Theta_{X_{\CC}/B_{\CC}}}^{-1}$ on $\TT^{-1}$. 
Then the norm of $u_L$ is independent of $L$. 
\end{Proposition}

\Proof
Let $b \in B(\CC)$. 
Since determinant line bundles are compatible with a base change 
and since the Quillen metric is given fiberwise, we get
\[
u_L : \det Rf_{b*}(L_b \otimes A_b) 
      \to \OO_{\Pic^{g-1}_{X_b}}(-\Theta_{X_b})\vert_{[L_b \otimes A_b]}, 
\]
where $[L_b\otimes A_b]$ is the point corresponding to $L_b \otimes A_b$ 
on $\Pic^{g-1}_{X_b}$. 
Then by the following lemma, 
we obtain Proposition~\ref{prop:key:lemma}. 
\QED

\begin{Lemma}
Let $M$ be a compact Riemann surface of genus $g \geq 1$, 
$L$ a line bundle of degree $0$ on $M$. 
We endow a $C^{\infty}$-metric $h_A$ on $A$,
a $C^{\infty}$-metric $h_{\Omega_M^1}$ on $\Omega_M^1$, 
and an admissible metric $h_L$ on $L$. 
Then we have a canonical isomorphism 
\[
u_L : \det \Gamma(L \otimes A) 
      \to \OO_{\Pic^{g-1}_{M}}(-\Theta_{M})\vert_{[L \otimes A]}, 
\]
where $\det \Gamma(L \otimes A)$ 
is the determinant line bundle of $L \otimes A$. 
We endow the Quillen metric on $\det \Gamma(L \otimes A)$ 
and $\Vert \cdot \Vert_{\Theta_M}^{-1}$ on  
$\OO_{\Pic^{g-1}_{M}}(-\Theta_{M})$. 
Then the norm of $u_L$ is independent of $L$. 
\end{Lemma}

\Proof
Let $h_A'$ and $h_{\Omega_M^1}'$ be admissible metrics on 
$A$ and $\Omega_M^1$ respectively. 
We write the Quillen metric defined by 
$(L \otimes A, h_L \otimes h_A)$ and 
$(\Omega_M^1, h_{\Omega_M^1})$ as 
$h_Q^{\overline{L}\otimes \overline{A}}$. 
We also write the Quillen metrics 
defined by $(L \otimes A, h_L \otimes h_A')$ and 
$(\Omega_M^1, h_{\Omega_M^1}')$ as 
$h_Q^{\overline{L}\otimes \overline{A}'}$. 
We decompose $u_L$ into 
\begin{multline*}
\left(\det \Gamma(L \otimes A), h_Q^{\overline{L}\otimes \overline{A}}\right) 
\overset{\alpha}{\longrightarrow}
\left(\det \Gamma(L \otimes A), h_Q^{\overline{L}\otimes \overline{A}'}\right) 
\\
\overset{\beta}{\longrightarrow}
\left(\det \Gamma(L \otimes A), h_F^{L\otimes A}\right) 
\overset{\gamma}{\longrightarrow}
\OO_{\Pic^{g-1}_{M}}(\Theta_{M})\vert_{[L \otimes A]}, 
\end{multline*}
where $h_F^{L\otimes A}$ is the Faltings' metric on 
$L \otimes A$. 
By the definition of the Quillen metrics, 
the norm of $\alpha$ is independent of $L$, 
because we only change the metric of $A$. 
The norm of $\beta$ is the difference of 
the Quillen metric and the Faltings' metric 
for admissible line bundles, which is a constant 
depending only on $M$ (cf. \cite[4.5]{So}). 
Moreover, the norm of 
$\gamma$ is also independent of $L$, 
which is actually given by $\exp\left( \delta(M)/8 \right)$ 
with the Faltings' delta function $\delta(M)$ 
(Or rather, this is the definition of $\delta(M)$) . 
Therefore the norm of $u_L$ 
is independent of $L$.  
\QED

%%%%%%%%%%%%%%%%%%%%%%%%%%%%%%%%%%%%%%%%%%%%%%%%%%%%%%%%
%%% Arithmetic height functions over function fields %%%
%%%%%%%%%%%%%%%%%%%%%%%%%%%%%%%%%%%%%%%%%%%%%%%%%%%%%%%%

\medskip
\section{Arithmetic height functions over function fields}
A. Moriwaki \cite{MoA} has recently constructed a theory of 
arithmetic height functions over function fields, 
with which he recovered the original Raynaud theorem 
(i.e., over a finitely generated field over $\QQ$).  
In this section, we see a part of his theory. 

Let $K$ be a finitely generated field over $\QQ$ with $\td_K(\QQ) =d$.  
Let $B$ be a normal projective arithmetic variety 
with the function field $K$. 
Let $\overline{H}$ be a {\em nef} $C^{\infty}$-hermitian $\QQ$-line bundle 
on $B$, i.e.,  
$\adeg(\overline{H}\vert_{C}) \geq 0$ for any curve $C$ 
and $\cherncl_1(\overline{H})$ is semi-positive on $B(\CC)$. 
A pair $\overline{B} = (B,\overline{H})$ with the above properties 
is called a {\em polarization} of $K$. 
Moreover, we say that a polarization $\overline{B}$ is {\em big} 
if $\rank H^0(B,H^{\otimes m})$ grows the order of $m^d$ and that 
there is a non-zero section $s$ of $H^0(B,H^{\otimes n})$  
with $\Vert s \Vert_{\sup} < 1$ for some positive integer $n$.  

Let $X_K$ be a projective variety over $K$
and $L_K$ a line bundle on $X_K$. 
By a {\em model} of $(X_K, L_K)$ over $B$, 
we mean a pair $(X \overset{f}{\longrightarrow} B, \overline{L})$ 
where $f: X \to B$ is a projective morphism of arithmetic varieties  
and $\overline{L} = (L, h_L)$ is a 
$C^{\infty}$-hermitian $\QQ$-line bundle on $X$ such that, 
on the generic fiber, $X$ and $L$ 
coincide with $X_K$ and $L_K$ respectively.  

By abbreviation, 
a model $(X \overset{f}{\longrightarrow} B, \overline{L})$ 
is sometimes written as $(X, \overline{L})$.
We note that although we use the notation $X_K$ and $L_K$, 
a model of $(X_K, L_K)$ is not a priori determined.  

For $P \in X(\overline{K})$, 
we denote by $\Delta_P$ the Zariski closure of
the $\Image\left( \Spec(\overline{K}) \to X_K \right)$ in $X$.  
Then we define the height of $P$ with respect to 
$(X \overset{f}{\longrightarrow} B, \overline{L})$ to be 
\[
h_{(X, \overline{L})}^{\overline{B}}(P)
= 
\frac{1}{[K(P):K]} 
\adeg\left( 
     \acherncl_1(\overline{L}\vert_{\Delta_P}) \cdot 
     \acherncl_1(f^*\overline{H}\vert_{\Delta_P})^d
     \right). 
\] 
If we change models of $(X_K,L_K)$, 
then height functions differ 
by only bounded functions on $X_K(\overline{K})$. 
Namely, 
if $(X, \overline{L})$ and $(X', \overline{L'})$ are 
two models of $(X_K,L_K)$, 
then there is a constant $C >0$ with 
\begin{equation}
\label{eqn:chande:of:model}
\left\vert h_{(X, \overline{L})}^{\overline{B}}(P) 
      -  h_{(X', \overline{L'})}^{\overline{B}}(P) \right\vert
\leq C
\end{equation}
for all $P \in X_K(\overline{K})$ (\cite[Corollary 3.3.5]{MoA}). 
Thus the height associated with $L_K$ and $\overline{B}$ 
is well-defined up to bounded functions on $X_K(\overline{K})$. 
We denote $h_{L_K}^{\overline{B}}$ 
the class of $h_{(X, \overline{L})}^{\overline{B}}$ 
modulo bounded functions.  

Now let $L_{\overline{K}}$ be a line bundle on 
$X_{\overline{K}} = X \otimes_K \overline{K}$. 
We would like to define 
$h_{L_{\overline{K}}}^{\overline{B}} : X_{\overline{K}} \to \RR$. 
For this, we need the following lemma. 

\begin{Proposition}
\label{prop:change:of:base}
Let $K'$ be a finite extension field of $K$, 
and let $g : B' \to B$ be a morphism of 
projective normal arithmetic varieties
such that the function field of $B'$ is $K'$. 
Let $X'$ be the main component of $X \times_B B'$ and 
\[
\begin{CD}
X' @>{g'}>> X \\
@V{f'}VV @V{f}VV \\
B' @>{g}>> B
\end{CD}
\]
the induced morphism. 
Then 
$h^{(B', g^*(\overline{H}))}_{X', g^*(\overline{L})}
= [K' : K] \  h^{(B, \overline{H})}_{(X,\overline{L})}$.   
\end{Proposition}

\Proof
\cite[Proposition 3.3.1]{MoA}
\QED

Let $L_{\overline{K}}$ be a line bundle on $X_{\overline{K}}$. 
We take a finite extension field $K'$ of $K$ such that 
$L_{\overline{K}}$ is defined over $X_{K'}$. 
Take a projective normal arithmetic variety $B'$ such that 
there is a morphism $g : B' \to B$ and that the function field 
of $B'$ is $K'$. Let $X'$ be the main component of $X \times_B B'$. 
We take a blow-up $\widetilde{X'} \to X'$ if necessary so that 
$L_{\overline{K}}$ extends to a line bundle $\widetilde{L'}$ on 
$\widetilde{X'}$. 

Then we define
\[
h^{\overline{B}}_ {L_{\overline{K}}}
=
\frac{1}{[K' : K]} \ 
h^{(B', g^*(\overline{H}))}_{(\widetilde{X'}, \widetilde{L'})}. 
\]
By \eqref{eqn:chande:of:model} 
and Proposition~\ref{prop:change:of:base}, 
it is easy to see that 
$h^{\overline{B}}_ {L_{\overline{K}}}$ is well-defined up to 
bounded functions on $X_{\overline{K}}(\overline{K})$. 
Moreover, if $L_{\overline{K}}$ is defined over $X_K$, 
then $h^{\overline{B}}_ {L_{\overline{K}}}$ is equal to 
$h^{\overline{B}}_{L_K}$. 

The next theorem shows some fundamental properties of 
$h_{L_{\overline{K}}}^{\overline{B}}$.  

\begin{Theorem}
\label{theorem:MoL:Northcott}
\begin{enumerate}
\renewcommand{\labelenumi}{(\roman{enumi})}
\item
(positiveness)
If we denote 
$\Supp \left( \Coker (H^0(X_{\overline{K}},L_{\overline{K}}) 
\otimes \OO_{X_{\overline{K}}}) \to L_{\overline{K}} \right)$ 
by $\Bs (L_{\overline{K}})$, then 
$h^{\overline{B}}_{L_{\overline{K}}}$ is bounded below on 
$\left(X_{\overline{K}} \setminus \Bs(L_{\overline{K}}) \right)$. 
\item
(Northcott)
Assume $\overline{H}$ is big and that 
$L_{\overline{K}}$ is ample.   
Then for any $e \geq 1$ and $M \geq 0$, 
\[
\{ P \in X_{\overline{K}}(\overline{K}) 
\mid  h^{\overline{B}}_{L_{\overline{K}}}(P) \leq M, \quad [K(P):K] \leq e \}
\] 
is a finite set. 
\end{enumerate}
\end{Theorem}

\Proof
c.f. \cite[Proposition 3.3.6 and Theorem 4.3]{MoA}
\QED

If $X_{\overline{K}}$ is an abelian variety, 
we can choose the good representative 
of a class $h_{L_{\overline{K}}}^{\overline{B}}$.  
For a line bundle $L_{\overline{K}}$ 
on $X_{\overline{K}}$ and a point $P \in X_{\overline{K}}(\overline{K})$, 
define $q_{L_{\overline{K}}}^{\overline{B}}(P,P)$ 
and $l_{L_{\overline{K}}}^{\overline{B}}(P)$ 
to be 
\begin{gather*}
q_{L_{\overline{K}}}^{\overline{B}}(P,P)
= \lim_{n \to \infty} \frac{1}{4^n} 
h^{\overline{B}}_{L_{\overline{K}}}(2^n P) \\
l_{L_{\overline{K}}}^{\overline{B}}(P) 
= \lim_{n \to \infty} \frac{1}{2^n} 
    \left(\frac{1}{4^n} h^{\overline{B}}_{L_{\overline{K}}}(2^n P) 
          - q_{L_{\overline{K}}}^{\overline{B}}(P,P)
    \right).
\end{gather*}
Then $q_{L_{\overline{K}}}^{\overline{B}}$ is a bilinear form, 
while $l_{L_{\overline{K}}}^{\overline{B}}$ is a linear form. 
We define $\widehat{h}^{\overline{B}}_{L_{\overline{K}}}$ by 
\[
\widehat{h}^{\overline{B}}_{L_{\overline{K}}}(P) 
= q_{L_{\overline{K}}}^{\overline{B}}(P,P) 
+ l_{L_{\overline{K}}}^{\overline{B}}(P),
\]
and call it the {\em canonical height} of $L_{\overline{K}}$ 
with respect to a polarization $\overline{B}$. 

\begin{Proposition}
\label{prop:MoL:canoanical:height}
Let $X_{\overline{K}}$ be an abelian variety. 
\begin{enumerate}
\renewcommand{\labelenumi}{(\roman{enumi})}
\item
If $L_{\overline{K}}$ is ample and symmetric, 
then $\widehat{h}^{\overline{B}}_{L_{\overline{K}}} \geq 0$. 
\item
If $L_{\overline{K}}$ and $M_{\overline{K}}$ are 
two line bundles on $X_{\overline{K}}$, 
then 
\[
\widehat{h}^{\overline{B}}_{L_{\overline{K}} \otimes M_{\overline{K}}}(P)
= \widehat{h}^{\overline{B}}_{L_{\overline{K}}}(P) 
+ \widehat{h}^{\overline{B}}_{M_{\overline{K}}}(P)
\]
\item
If $P$ is a torsion point, then 
$\widehat{h}^{\overline{B}}_{L_{\overline{K}}}(P) = 0$. 
If we assume $\overline{H}$ is big, 
then $\widehat{h}^{\overline{B}}_{L_{\overline{K}}}(P) = 0$ 
if and only if $P$ is a torsion point. 
\end{enumerate}
\end{Proposition}

\Proof
The first assertion follows from Theorem~\ref{theorem:MoL:Northcott}(i). 
The second assertion can be readily checked. 
The third assertion is an easy consequence of 
Theorem~\ref{theorem:MoL:Northcott}(ii). 
We note that in (i) we need the symmetricity 
of a line bundle. 
\QED

We need the next lemma to prove Proposition~\ref{prop:main}. 

\begin{Lemma}
\label{lemma:thm:of:square}
Let $L_{\overline{K}}$ is an ample symmetric line bundle on an 
abelian variety $X_{\overline{K}}$, 
P an element of $X_{\overline{K}}(\overline{K})$. 
Let $t$ be an element of $X_{\overline{K}}(\overline{K})$ 
and $T_t : X_{\overline{K}} \to X_{\overline{K}}$ the translation by $t$.   
Then there is a constant $C$ such that 
\[
\left\vert
\widehat{h}^{\overline{B}}_{T_t^*(L_{\overline{K}})}(n P)
- n^2 \widehat{h}^{\overline{B}}_{L_{\overline{K}}}(P)
\right\vert
= C n
\] 
for any positive integers n. 
\end{Lemma}

\Proof
Let $T_{-t} : X_{\overline{K}} \to X_{\overline{K}}$ 
be the translation by $-t$. 
We write $T_t^*(L_{\overline{K}})^{\otimes 2}$ as 
\[
T_t^*(L_{\overline{K}})^{\otimes 2} 
= \left( T_t^*(L_{\overline{K}}) \otimes T_{-t}^*(L_{\overline{K}}) \right)
  \otimes \left( T_t^*(L_{\overline{K}}) 
  \otimes (T_{-t}^*(L_{\overline{K}}))^{-1} \right)
\]
Since $T_t^*(L_{\overline{K}}) \otimes T_{-t}^*(L_{\overline{K}}) 
= L_{\overline{K}} ^{\otimes 2}$ 
by the theorem of square, 
we obtain 
\[
T_t^*(L_{\overline{K}})^{\otimes 2} 
= \left( L_{\overline{K}}^{\otimes 2} \right)
  \otimes \left( T_t^*(L_{\overline{K}}) 
  \otimes (T_{-t}^*(L_{\overline{K}}))^{-1} \right). 
\]
Thus we get 
$4 \widehat{h}^{\overline{B}}_{T_t^*(L_{\overline{K}})} 
= 4 \widehat{h}^{\overline{B}}_{L_{\overline{K}}} + 
\widehat{h}^{\overline{B}}_{ T_t^*(L_{\overline{K}}) 
\otimes (T_{-t}^*(L_{\overline{K}}))^{-1}}$. 
Since $L_{\overline{K}}$ is symmetric 
and $T_t^*(L_{\overline{K}}) \otimes (T_{-t}^*(L_{\overline{K}}))^{-1}$
is anti-symmetric, 
$\widehat{h}^{\overline{B}}_{L_{\overline{K}}}$ is quadric, 
while $\widehat{h}^{\overline{B}}_{ T_t^*(L_{\overline{K}}) 
\otimes (T_{-t}^*(L_{\overline{K}}))^{-1}}$ 
is linear. 
Thus if we set 
$C = \left\vert 
\widehat{h}^{\overline{B}}_{ T_t^*(L_{\overline{K}}) 
\otimes (T_{-t}^*(L_{\overline{K}}))^{-1}}(P) 
\right\vert$, then  
we obtain the lemma. 
\QED

%%%%%%%%%%%%%%%%%%%%%%%%%%%%%%%
%%% Height and intersection %%%
%%%%%%%%%%%%%%%%%%%%%%%%%%%%%%%

\medskip
\section{Height and intersection}
By a {\em big} Zariski open set of a noetherian scheme $B$, 
we mean a Zariski open set $B'$ of $B$ 
with $\codim_{B}(B\setminus B') \geq 2$. 

We first prove the following proposition, 
which is a special case of the main theorem 
(Theorem~\ref{theorem:main}). 

\begin{Proposition}
\label{prop:main}
Let $K$ be a finitely generated field over $\QQ$, 
$X_K$ a geometrically irreducible regular projective curve over $K$, 
and $L_K$ a line bundle on $X_K$ with $\deg L_K =0$. 
Let $\overline{B} = (B, \overline{H})$ be a polarization of $K$, 
and $(X \overset{f}{\longrightarrow} B, \overline{L})$ 
a model of $(X_K, L_K)$. 
We make the following assumptions on the model:
\begin{enumerate}
\renewcommand{\labelenumi}{(\alph{enumi})}
\item
$B$ is regular;
\item
$f$ is semi-stable with a section $\epsilon$;  
\item
$X_{\CC}$ and $B_{\CC}$ are non-singular and 
$f_{\CC} : X_{\CC} \to B_{\CC}$ is smooth.  
\end{enumerate}
Let $J_K$ be the Jacobian of $X_K$ and 
$\Theta_{\overline{K}}$ a divisor on $J_{\overline{K}}$ 
which is a translation of 
the theta divisor on $\Pic^{g-1}(X_{\overline{K}})$ 
by a theta characteristic.  
If there is a big Zariski open set $B' \subset B$ 
such that $\deg(L\vert_C)=0$ for any fibral curve $C$ 
lying over $B'$ and 
if the metric of $\overline{L}$ is flat along fibers, 
then 
\begin{equation}
\label{eqn:main:1}
\adeg\left(\acherncl_1(\overline{L})^2 \cdot 
           \acherncl_1(f^*(\overline{H}))^d 
     \right)
= 
-2 \widehat{h}_{\OO_{J_{\overline{K}}}(\Theta_{\overline{K}})}^{\overline{B}}([L_K]), 
\end{equation}
where $[L_K]$ denotes the point of $J_K$ corresponding to $L_K$. 
\end{Proposition}

\Proof
We note that since $\deg L = 0$, 
the admissibility of $\overline{L}$ means that 
the metric of $\overline{L}$ is flat along fibers. 
Since $\deg(L_K) = 0$, if we change $L$ to $L \otimes f^*(M)$ 
with $M$ being a line bundle on $B$, 
then the right hand of \eqref{eqn:main:1} does not change. 
Thus we may assume that $L$ is rigidified along the section $\epsilon$.  
Let us set $A = \OO_X((g-1)[\epsilon])$. 
Then $A$ is a rigidified line bundle  of degree $(g-1)$ on $X$. 
Let $(P^a, U^a)$ be the translation of 
$(\Pic^0_{X/B}, U^{0})$ by $A$, 
where $U^{0}$ is the universal line bundle on 
$X \times_{B} \Pic^0_{X/B}$. 
We put $\TT^{-1} = \det Rq^a_*(U^a)$, 
where $q^a :  X \times_{B} P^a \to P^a$ 
is the second projection. 

We give an admissible metric $h_A$ on $A$ 
and an admissible metric $h_{\omega_{X/B}}$ on $\omega_{X/B}$ 
and then give $\det Rf_*(L^{\otimes n} \otimes A)$ 
the Quillen metric 
$h_Q^{\overline{L}^{\otimes n} \otimes \overline{A}}$ 
with respect to $\overline{L}^{\otimes n} \otimes \overline{A}
= (L^{\otimes n} \otimes A, h_L^{n} \cdot h_A)$ 
and $\overline{\omega_{X/B}} = (\omega_{X/B},h_{\omega_{X/B}})$. 
Moreover we endow $\Vert \cdot \Vert_{\Theta_{X_{\CC}/B_{\CC}}}^{-1}$ 
on $\TT^{-1}$ (cf. Proposition~\ref{prop:metric:of:theta}). 

Let us put $X' = f^{-1}(B')$, $f' = f\vert_{X'}$, 
$L' = L\vert_{X'}$ and $A' = A\vert_{X'}$. 
Moreover Let 
$(P^a{}', U^a{}')$, $(\Pic^0_{X'/B'}, U^{0}{}')$, 
$q^a{}'$ and $\TT^{-1}{}' = \det Rq^a{}'_*(U^a{}')$  
be the restriction of 
$(P^a, U^a)$, $(\Pic^0_{X/B}, U^{0})$, 
$q^a$ and $\TT^{-1}{} = \det Rq^a_*(U^a{})$ over $B'$, respectively. 

Now we consider $L'{}^{\otimes n} \otimes A'$ for a positive integer $n$. 
Since $\deg(L'\vert_{C}) = 0$ 
for any fibral curve lying over $B'$, 
$L'$ belongs to $\Pic^0_{X'/B'}$. 
Thus by Theorem~\ref{theorem:det:for:semistable}(ii), 
there is a canonical morphism 
$g_n' : B' \to P^a{}'$ such that 
\[
u_n' : \det Rf_*'(L'{}^{\otimes n} \otimes A') 
       \overset{\sim}{\longrightarrow} g_n'{}^*(\TT^{-1}{}')
\]
is canonically isomorphic over $B'$. 
Since both sides are metrized, 
we can consider the norm $\alpha_n$ of $u_n'$. 
Then
\[
u_n' : \left(\det Rf_*'(L'{}^{\otimes n} \otimes A'),
             h_Q^{\overline{L}^{\otimes n} \otimes \overline{A}}
       \right) 
       \overset{\sim}{\longrightarrow} 
       g_n'{}^*\left(\TT^{-1}{}',
       \Vert \cdot \Vert_{\Theta_{X_{\CC}/B_{\CC}}}^{-1}
       \right) 
       \otimes \OO_{B'}\left( \alpha_n^{-1} \right)
\]
is an isometry. 
Moreover, by Proposition~\ref{prop:key:lemma}, 
the function $\alpha_n : B_{\CC}(\CC) \to \RR_{> 0}$ 
is independent of $n$. 

Next we consider a compactification of $P^a$. 
Since there is a relatively ample line bundle on $P^a$, 
we first embed $P^a$ into a large projective space $\PP_B^N$ 
and then take its closure. 
If $\TT^{-1}$ does not extend to a line bundle on this closure, 
then we make blow-ups along the boundary. 
Then we get a projective arithmetic variety $\underline{P}^a$ 
with $\pi : \underline{P}^a \to B$ and 
a line bundle $\underline{\TT}^{-1}$ on $\underline{P}^a$ 
with $\underline{\TT}^{-1}\vert_{P^a} = \TT^{-1}$. 
We note that since $f_{\CC}$ is smooth, 
$\underline{P}^a_{\CC} = P^a_{\CC}$ 

Let $\Delta_n$ be the Zariski closure of the 
$\Image(g_n' : B' \to P^a{}')$ in $\underline{P}^a$. 
Now we claim the following equation; 
\begin{multline}
\label{eqn:main:two}
\adeg\left( 
         \acherncl_1\left(\det Rf_*(L^{\otimes n} \otimes A),
             h_Q^{\overline{L}^{\otimes n} \otimes \overline{A}}\right)
         \cdot \acherncl_1(\overline{H})^d    
      \right) \\
= \adeg\left( 
     \acherncl_1
     \left(
     \OO_{\underline{P}^a}(\underline{\TT}^{-1}), 
     \Vert \cdot \Vert_{\Theta_{X_{\CC}/B_{\CC}}}^{-1}
     \right)\left\vert_{\Delta_n}\right. 
     \cdot \acherncl_1(\pi^*(\overline{H}))^d\left\vert_{\Delta_n}\right.
     \right)
-       
\frac{1}{2}\int_{B_{\CC}(\CC)} (\log \alpha_n) 
\wedge \cherncl_1(\overline{H})^d. 
\end{multline}

Actually, since $B$ is regular and $B'$ is big, 
a line bundle on $B'$ extends uniquely to a line bundle on $B$. 
The line bundle $\det Rf_*'(L'{}^{\otimes n} \otimes A')$ on $B'$ 
extends to $\det Rf_*(L^{\otimes n} \otimes A)$ 
and the line bundle $g_n'{}^*(\TT^{-1}{}')$ on $B'$ 
extends to a line bundle on $B$, which we denote by 
$M_n$. Let us set 
$\overline{M_n} 
= (M_n, g_n'{}^*(\Vert \cdot \Vert_{\Theta_{X_{\CC}/B_{\CC}}}^{-1}))$. 
Since $\pi \vert_{\Delta_n} : \Delta_n \to B$ is an isomorphism 
over $B'$ and $\codim_{B}(B \setminus B') \geq 2$, 
$\overline{M_n}$ is actually equal to 
$(\pi \vert_{\Delta_n})_*(\OO_{\underline{P}^a}(\underline{\TT}^{-1}), 
     \Vert \cdot \Vert_{\Theta_{X_{\CC}/B_{\CC}}}^{-1})$. 
Then since the infinite part is not altered at all, 
we get the isometry
\[
u_n : \left(\det Rf_*(L^{\otimes n} \otimes A),
             h_Q^{\overline{L}^{\otimes n} \otimes \overline{A}}
       \right) 
       \overset{\sim}{\longrightarrow} 
       \overline{M_n} 
       \otimes \OO_{B}\left( \alpha_n^{-1} \right).
\]
Then by intersecting $\acherncl_1(\overline{H})^d$ 
and taking degrees on both sides, we get
\begin{multline*}
\adeg\left( 
         \acherncl_1\left(\det Rf_*(L^{\otimes n} \otimes A),
             h_Q^{\overline{L}^{\otimes n} \otimes \overline{A}}\right)
         \cdot \acherncl_1(\overline{H})^d    
      \right) \\
= 
\adeg\left(        
       \acherncl_1(\overline{M_n})
       \cdot \acherncl_1(\overline{H})^d
     \right) 
-       
\frac{1}{2}\int_{B_{\CC}(\CC)} (\log \alpha_n) 
\wedge \cherncl_1(\overline{H})^d  \\
= 
\adeg\left( 
     \acherncl_1
     \left(
     \OO_{\underline{P}^a}(\underline{\TT}^{-1}), 
     \Vert \cdot \Vert_{\Theta_{X_{\CC}/B_{\CC}}}^{-1}
     \right)\vert_{\Delta_n} 
     \cdot \acherncl_1(\pi^*(\overline{H}))^d\vert_{\Delta_n}
     \right)
-       
\frac{1}{2}\int_{B_{\CC}(\CC)} (\log \alpha_n) 
\wedge \cherncl_1(\overline{H})^d, 
\end{multline*}
where we use the projection formula in the second equality. 

First we compute the left hand side of \eqref{eqn:main:two}. 
By the arithmetic Riemann-Roch theorem 
established by Gill\'{e}t and Soul\'{e} \cite{GS}, 
we have 
\begin{multline*}
\acherncl_1\left(\det Rf_*(L^{\otimes n} \otimes A),
             h_Q^{\overline{L}^{\otimes n} \otimes \overline{A}}\right) \\
= \frac{1}{2} f_* \left( 
  \acherncl_1(\overline{L}^{\otimes n} \otimes \overline{A})^2 
 -\acherncl_1(\overline{L}^{\otimes n} \otimes \overline{A}) \cdot 
  \acherncl_1(\overline{\omega_{X/B}})       
                  \right)
+ \acherncl_1\left(
      \det Rf_*(\OO_X), h_Q^{\overline{\OO_X}}
             \right) \\
= \frac{1}{2} f_* \left( \acherncl_1(\overline{L})^2 \right) n^2 + O(n). 
\end{multline*}

Thus, we obtain
\begin{multline}
\label{eqn:main:left}
\adeg\left( 
         \acherncl_1\left( 
                  \det Rf_*(L^{\otimes n} \otimes A),
                  h_Q^{\overline{L}^{\otimes n} \otimes \overline{A}}
                    \right)
         \cdot \acherncl_1(\overline{H})^d    
      \right) \\
=  \frac{1}{2} \adeg 
   \left(
   f_* \left( \acherncl_1(\overline{L})^2 \right)
   \cdot \acherncl_1(\overline{H})^d   
   \right) n^2 + O(n) \\
= \frac{1}{2} \adeg 
   \left(
   \acherncl_1(\overline{L})^2 
   \cdot \acherncl_1(f^*(\overline{H}))^d   
   \right) n^2 + O(n) 
\end{multline}

Next we compute the right hand side of \eqref{eqn:main:two}. 
Let $\lambda_a : \Pic_{X/B}^{0} 
\overset{\sim}{\longrightarrow} P^a$ be the isomorphism 
which is given by 
the translation by $A$. 
By way of this identification, 
let $\underline{P}^0$ be the compactification of 
$\Pic_{X/B}^0$ which corresponds to $\underline{P}^a$. 
Similarly, we define   
$(\underline{\TT}^0)^{-1}$, $\Delta_n^0$ and $\pi^0$ 
which correspond to 
$\underline{\TT}^{-1}$, $\Delta_n$ and $\pi$ 
respectively. 
We note that a metric on $(\underline{\TT}^0)^{-1}$ 
induced from $\lambda_a$ is nothing but 
$\Vert \cdot \Vert_{\Theta^0_{X_{\CC}/B_{\CC}}}^{-1}$ 
by Proposition~\ref{prop:metric:of:theta}. 
Then we have  
$\underline{\TT}^0_K = \OO_{J_K}(\Theta_K')$,  
where 
\[
\Theta_K' 
= \Theta_K + [\text{a theta characteristic}] - (g-1)[\epsilon_K].  
\] 

Since $\left(\pi^0 : \underline{P}^0 \to B,
\left( (\underline{\TT}^0)^{-1},
\Vert \cdot \Vert_{\Theta_{X_{\CC}/B_{\CC}}^0}^{-1} \right) \right)$ 
is a model of $\left(J_K, \OO_{J_K}(-\Theta_K')\right)$, 
\eqref{eqn:chande:of:model} shows that 
there is a constant $C$ such that 
\begin{multline*}
\left\vert 
\adeg\left( 
     \acherncl_1
     \left(
     \OO_{\underline{P}^a}(\underline{\TT}^{-1}), 
     \Vert \cdot \Vert_{\Theta_{X_{\CC}/B_{\CC}}}^{-1}
     \right)\left\vert_{\Delta_n}\right. 
     \cdot \acherncl_1(\pi^*(\overline{H}))^d 
     \left\vert_{\Delta_n}\right.
     \right)
- 
\widehat{h}_{\OO_{J_K}(-\Theta_K')}^{\overline{B}}([L_K^{\otimes n}])
\right\vert \\
= 
\left\vert
\adeg\left( 
     \acherncl_1
     \left(
     \OO_{\underline{P}^0}((\underline{\TT}^0)^{-1}), 
     \Vert \cdot \Vert_{\Theta_{X_{\CC}/B_{\CC}}^0}^{-1}
     \right)\left\vert_{\Delta_n}\right. 
     \cdot \acherncl_1(\pi^*(\overline{H}))^d
     \vert_{\Delta_n}
     \right)
- 
\widehat{h}_{\OO_{J_K}(-\Theta_K')}^{\overline{B}}([L_K^{\otimes n}])
\right\vert 
\leq C 
\end{multline*}
Then using Lemma~\ref{lemma:thm:of:square}, 
we get 
\begin{equation}
\label{eqn:main:right} 
\left\vert 
\adeg\left( 
     \acherncl_1
     \left(
     \OO_{\underline{P}^a}(\underline{\TT}^{-1}), 
     \Vert \cdot \Vert_{\Theta_{X_{\CC}/B_{\CC}}}^{-1}
     \right)\left\vert_{\Delta_n}\right. 
     \cdot \acherncl_1(\pi^*(\overline{H}))^d
     \vert_{\Delta_n}
     \right)
- 
n ^2 \widehat{h}_{\OO_{J_{\overline{K}}}
(\Theta_{\overline{K}})}^{\overline{B}}([L_K])
\right\vert
= O(n). 
\end{equation}

Taking into consideration 
\eqref{eqn:main:left} and \eqref{eqn:main:right} 
and the fact that $\alpha_n$ is independent of $n$, 
if we divide \eqref{eqn:main:two} by $n^2$ 
and let $n$ goes to $\infty$, 
we get \eqref{eqn:main:1}. 
\QED

Now we prove the main theorem of this paper. 

\begin{Theorem}
\label{theorem:main}
Let $K$ be a finitely generated field over $\QQ$, 
$X_K$ a geometrically irreducible regular projective curve over $K$, 
and $L_K$ a line bundle on $X_K$ with $\deg L_K =0$. 
Let $\overline{B} = (B, \overline{H})$ be a polarization of $K$, 
and $(X \overset{f}{\longrightarrow} B, \overline{L})$ 
a model of $(X_K, L_K)$. 
We make the following assumptions on the model:
\begin{enumerate}
\renewcommand{\labelenumi}{(\alph{enumi})}
\item
$f$ is semi-stable;  
\item
$X_{\CC}$ and $B_{\CC}$ are non-singular and 
$f_{\CC} : X_{\CC} \to B_{\CC}$ is smooth.  
\end{enumerate}

Then we have 
\[
\adeg\left( \acherncl_1(\overline{L})^2 \cdot 
            \acherncl_1\left(f^*(\overline{H})\right)^{d}
     \right)
\leq 
-2 \widehat{h}_{\OO_{J_{\overline{K}}}
(\Theta_{\overline{K}})}^{\overline{B}}([L_K]), 
\]
where $[L_K]$ denotes the point of $J_K$ corresponding to $L_K$. 

Furthermore, we assume 
that $H$ is ample and $\cherncl_1(\overline{H})$ is positive. 
Then the equality holds 
if and only if $\overline{L}$ satisfies the following properties:
\begin{enumerate}
\renewcommand{\labelenumi}{(\alph{enumi})}
\item 
There is a big Zariski open set $B''$ of $B$ 
such that $\deg(L\vert_C) = 0$ for any fibral curves $C$ lying over $B''$.
\item 
The metric of $\overline{L}$ is flat along fibers. 
\end{enumerate}
\end{Theorem}

The next corollary is an immediate consequence of 
the main theorem and Proposition~\ref{prop:MoL:canoanical:height}(iii). 

\begin{Corollary}
\label{cor:of:main:thm}
Let the notation and the assumption be as in Theorem~\ref{theorem:main}. 
We assume that $\overline{H}$ is big, $H$ is ample and 
$\cherncl_1(H)$ is positive. 
Then 
\[
\adeg\left( \acherncl_1(\overline{L})^2 \cdot 
            \acherncl_1\left(f^*(\overline{H})\right)^{d}
     \right)
= 0
\]
if and only if the following properties hold:
\begin{enumerate}
\renewcommand{\labelenumi}{(\alph{enumi})}
\item 
There is a big Zariski open set $B''$ of $B$ 
such that $\deg(L\vert_C) = 0$ for any fibral curves $C$ lying over $B''$; 
\item 
The restriction of the metric of $\overline{L}$ 
to each fiber is flat; 
\item
There is a positive integer $m$ with 
$L_K^{\otimes m} = \OO_{X_K}$.  
\end{enumerate}

\end{Corollary}

We need three lemmas to prove the theorem. 

\begin{Lemma}
\label{lemma:for:main:reduction}
Let $\widetilde{K}$ be a finite extension field of $K$, 
and let $g : \widetilde{B} \to B$ be a morphism 
of projective normal arithmetic varieties 
such that the function field of $\widetilde{B}$ 
is $\widetilde{K}$. 
Let $\widetilde{X} = X \times_B \widetilde{B}$ 
and 
\[
\begin{CD}
\widetilde{X} @>{\widetilde{g}}>> X \\
@V{\widetilde{f}}VV @V{f}VV \\
\widetilde{B} @>{g}>> B
\end{CD}
\]
the induced morphism. 
Then 
\[
\adeg\left( \acherncl_1(\widetilde{g}^* \overline{L})^2 \cdot 
            \acherncl_1\left(\widetilde{f}^* g^* (\overline{H})\right)^{d}
     \right) 
= 
[\widetilde{K} : K] \ 
\adeg\left( \acherncl_1(\overline{L})^2 \cdot 
            \acherncl_1\left(f^*(\overline{H})\right)^{d}
     \right).
\]
\end{Lemma}

\Proof
It is an easy consequence of the projection formula. 
\QED

\begin{Lemma}
\label{lemma:for:main:infinite}
Let $\overline{L} = (L, h_L)$ 
be a $C^{\infty}$-hermitian line bundle on $X$ 
and $\overline{L}' = (L, h_L')$ 
be a hermitian line bundle whose metric is flat 
along fibers. Then 
\[
\adeg\left( \acherncl_1(\overline{L})^2 \cdot 
            \acherncl_1\left(f^*(\overline{H})\right)^{d}
     \right)
\leq
\adeg\left( \acherncl_1(\overline{L}')^2 \cdot 
            \acherncl_1\left(f^*(\overline{H})\right)^{d}
     \right)
\]
If $\cherncl_1(\overline{H})$ is positive 
over a dense open subset of $B(\CC)$, 
then the equality holds if and only if 
the metric of $\overline{L}$ is flat along fibers. 
\end{Lemma}

\Proof
Let us write $h_{L'} = u h_L$. 
Then $u$ is a positive smooth function on $X_{\CC}(\CC)$. Since
\[
\acherncl_1(\overline{L}) 
= \acherncl_1(\overline{L}') + (0, \log u),
\]
we have 
\[
\acherncl_1(\overline{L})^2
= \acherncl_1(\overline{L}')^2 
+ (0, 2 \cherncl_1(\overline{L}') \log u)
+ \left(0, (\log u) d d^c (\log u) \right).
\]
Thus 
\begin{multline*}
\adeg\left( \acherncl_1(\overline{L})^2 \cdot 
            \acherncl_1\left(f^*(\overline{H})\right)^{d}
     \right) 
=
\adeg\left( \acherncl_1(\overline{L}')^2 \cdot 
            \acherncl_1\left(f^*(\overline{H})\right)^{d}
     \right) \\
- \int_{X_{\CC}(\CC)} (\log u) \cherncl_1(\overline{L}') \wedge 
\cherncl_1\left(f^*(\overline{H})\right)^d
+ \frac{1}{2}\int_{X_{\CC}(\CC)} (\log u) d d^c (\log u) \wedge 
\cherncl_1\left(f^*(\overline{H})\right)^d.
\end{multline*}
Now the assertion follows from the following two claims.

\begin{Claim}
$\int_{X_{\CC}(\CC)} (\log u) \cherncl_1(\overline{L}') \wedge 
\cherncl_1\left(f^*(\overline{H})\right)^d = 0$
\end{Claim}

\Proof
For $b \in B_{\CC}(\CC)$, 
$ \cherncl_1(\overline{L}')\vert_b = 0$. Then 
\[
\int_{X_{\CC}(\CC)} (\log u) \cherncl_1(\overline{L}') \wedge 
\cherncl_1\left(f^*(\overline{H})\right)^d
=
\int_{B_{\CC}(\CC)} 
\left( \int_{f_{\CC} : X_{\CC} \to B_{\CC}} (\log u)  
\cherncl_1(\overline{L}') \right)
\cherncl_1\left(\overline{H}\right)^d = 0.
\]
\QED

\begin{Claim}
$\int_{X_{\CC}(\CC)} (\log u) d d^c (\log u) \wedge 
\cherncl_1\left(f^*(\overline{H})\right)^d \leq 0$. 
Moreover, if $\cherncl_1(\overline{H})$ is positive 
over a dense open set of $B(\CC)$, 
then the equality holds if and only if $u = f^*(v)$ 
with some $C^{\infty}$ function $v$ on $B_{\CC}(\CC)$. 
\end{Claim}

\Proof
We have 
\begin{align*}
(\log u) d d^c (\log u) 
 & = \frac{\sqrt{-1}}{2 \pi} (\log u) 
     \partial \overline{\partial} (\log u) \\
 & = \frac{\sqrt{-1}}{2 \pi}
     \partial \left( \log u \cdot \overline{\partial} (\log u) 
              \right)
    - \frac{\sqrt{-1}}{2 \pi} \partial ( \log u ) \wedge 
            \overline{\partial} (\log u).
\end{align*}
Since $\cherncl_1\left(f^*(\overline{H})\right)^d$ 
is a closed $(d,d)$-form, by Stokes' lemma, we get 
\begin{multline*}
\int_{X_{\CC}(\CC)} (\log u) d d^c (\log u) \wedge 
\cherncl_1\left(f^*(\overline{H})\right)^d \\
= 
-\frac{1}{2 \pi}
\int_{X_{\CC}(\CC)}  \left(\sqrt{-1} \partial ( \log u ) \wedge 
          \overline{\partial} (\log u) \right) \wedge 
\cherncl_1\left(f^*(\overline{H})\right)^d.
\end{multline*}
By the definition of the polarization 
of $\overline{B} = (B, \overline{H})$, $\cherncl_1(\overline{H})$ 
is semipositive. 
Moreover, $\partial ( \log u ) \wedge \overline{\partial} (\log u)$ 
is semipositive. Thus we get the first assertion.

Suppose now $\cherncl_1(\overline{H})$ is positive
over a dense open set of $B(\CC)$. We have
\begin{multline*}
\int_{X_{\CC}(\CC)}  \left(\sqrt{-1} \partial ( \log u ) \wedge 
          \overline{\partial} (\log u) \right) \wedge 
\cherncl_1\left(f^*(\overline{H})\right)^d \\
= 
\int_{B_{\CC}(\CC)}  
\left(
 \int_{f_{\CC} : X_{\CC} \to B_{\CC}} 
 \sqrt{-1} \partial ( \log u ) \wedge  
\overline{\partial} (\log u) \right) 
\cherncl_1\left(\overline{H}\right)^d.
\end{multline*}
If this value is zero, then, for any $b \in B_{\CC}$, 
$\sqrt{-1} \partial ( \log u ) \wedge  
\overline{\partial} (\log u) \vert_{X_b} =0$. 
Then $u\vert_{X_b}$ is a constant function on $X_b(\CC)$. 
This shows the second assertion.
\QED

\begin{Lemma}
\label{lemma:for:main:finite}
We assume that $B$ is regular. 
Let $\Delta$ be the set of critical values of $f$, i.e., $\Delta = 
\{ b \in B \mid \text{$f$ is not smooth over $b$} \}$. 
Let $\Delta = \bigcup_{i=1}^{I} \Delta_i$ be the irreducible decomposition 
of $\Delta$ such that $\Delta_1, \ldots, \Delta_{I_1}$ 
are divisors on $B$ while $\codim_B (\Delta_i) \geq 2$ for $i \geq I_1+1$. 
Let us set  $\Gamma_i = f^{-1} (\Delta_i)$ 
for $i = 1, \ldots, I_1$ and 
write $\Gamma_i = \bigcup_{j=1}^{J_i} \Gamma_{ij}$ 
as its irreducible decomposition. 
Note that $\Gamma_{ij}$ are all divisors on $X$ 
for $1 \leq i \leq I_1, 1 \leq j \leq J_i$. 
Then there are a big Zariski open set $B'$ of $B$, 
integers $e_{ij} \ (1 \leq i \leq I_1, 1 \leq j \leq J_i)$ 
and a positive integer $m$ 
such that $L^{\otimes m} \otimes \OO_X (- \sum_{ij} e_{ij} 
\Gamma_{ij}) \vert_{B'}$ belongs to $\Pic^0_{f^{-1}(B')/B'}$.
\end{Lemma}

\Proof
If $I_1 =0$, then we have nothing to prove. 
Thus, we assume $I_1 \geq 1$. 
To ease the notation, 
we first assume the irreducibility of $\Delta$. 
Since $f_{\CC}$ is smooth, 
$\Delta$ is defined over the finite field $\FF_p$ 
for some prime number $p$. 
Let $k(\Delta)$ be the rational function of $\Delta$ 
and write $\eta = \Spec(k(\Delta))$. 
Moreover, let $\overline{k(\Delta)}$ be an 
algebraic closure of $k(\Delta)$ and write 
$\overline{\eta} = \Spec(\overline{k(\Delta)})$. 

Let $X_{\overline{\eta}} = \cup_{ 1 \leq j \leq J} \cup_{1 \leq \alpha \leq \alpha(j)} 
C_{j}^{\alpha}$ be the irreducible decomposition of 
$X_{\overline{\eta}}$ such that $C_{j}^{\alpha}$ and $C_{j}^{\beta}$ are 
$\Gal(\overline{k(\Delta)}/k(\Delta))$-conjugate to each other 
for $1 \leq \alpha, \beta \leq \alpha(j)$. 
We denote by $\Gamma_j$ the Zariski closure of 
$C_{j}^{\alpha}$ in $X$ for some (hence all) $\alpha$. 

We put $c_{j}^{\alpha} = \deg (L_{\eta}\vert_{C_{j}^{\alpha}})$. 
Since $L$ is defined over $X$, 
$c_{j}^{\alpha} = c_{j}^{\beta}$ 
for $1 \leq \alpha, \beta \leq \alpha(j)$. 
Moreover, since the degree of $L$ is zero, 
$\sum_{1 \leq j \leq J, 1 \leq \alpha \leq \alpha(j)} c_{j}^{\alpha} = 0$. 

We put $q_{jk}^{\alpha \beta} 
= \dim_{k(\Delta)}(C_{j}^{\alpha} \cap C_{k}^{\beta})$ 
for $(j, \alpha) \neq (k, \beta)$, 
and $q_{jj}^{\alpha \alpha} = 
- \sum_{(k, \beta) \neq (j, \alpha)} q_{jk}^{\alpha \beta}$. 
Then by Zariski's lemma (\cite[I, Lemma~(2.10)]{Barth}), 
there are rational numbers 
$a_{j}^{\alpha} \ (1 \leq j \leq J, 1 \leq \alpha \leq \alpha(j))$ 
such that $a_{j}^{\alpha} = a_{j}^{\beta}$ and that 
$\sum_{j ,\alpha} a_{j}^{\alpha} q_{jk}^{\alpha \beta} 
= c_k^{\beta}$ for $1 \leq k \leq J$ 
and $1 \leq \beta \leq \alpha(k)$. 
Moreover, $\sum_{j, k, \alpha, \beta} 
a_{j}^{\alpha} a_ k^{\beta} q_{jk}^{\alpha \beta} = 0$ 
if and only if $a_{j}^{\alpha} = a_ k^{\beta}$ 
for any $(j, \alpha)$ and $(k, \beta)$. 

Let $Y$ be the subset of 
$\vert \Delta \vert$ 
consisting of $\overline{\FF_p}$-valued points $b$ such that:
\begin{enumerate}
\renewcommand{\labelenumi}{(\alph{enumi})}
\item
The irreducible decomposition of $X_b$ 
is of form 
$X_b = \cup_{ 1 \leq j \leq J} \cup_{1 \leq \alpha \leq \alpha(j)} 
C(b)_{j}^{\alpha}$ such that 
$\Gamma_j \cap X_b = \cup_{1 \leq \alpha \leq \alpha(j)} C(b)_{j}^{\alpha}$; 
\item
$\deg (L \vert_{C(b)_{j}^{\alpha}}) = c_{j}^{\alpha}$;
\item
$\Gamma_j \cdot C(b)_{k}^{\beta} 
= \sum_{1 \leq \alpha \leq \alpha(j)} q_{jk}^{\alpha \beta}$. 
\end{enumerate}
Then there is a divisor $Z$ on $\Delta$ 
such that $Y \subset \vert Z \vert$. We set $B' = B - \vert Z \vert$. 

Now we set $e_{j} = m a_{j}^{\alpha} \ (1 \leq j \leq J)$ 
for sufficiently divisible $m$ and 
$L' = L^{\otimes m} \otimes \OO_{X}(-\sum_{j=1}^{J}e_j\Gamma_j)$. 
We claim that $L' \vert_{B'}$ belongs to $\Pic^0_{f^{-1}(B')/B'}$. 
Indeed, if $b \not\in \Delta$, 
then $X_b$ is a smooth connected curve and $\deg(L'\vert_{X_b}) =0$. 
Thus $L'\vert_{X_b}$ belongs to $\Pic^0_{X_b}$. 
Next, if $b \in \Delta \setminus \vert Z \vert$, 
then $X_b = \cup_{j, \alpha} C(b)_{j}^{\alpha}$ 
is the irreducible decomposition of $X_b$ 
and 
\[
\deg (L'\vert_{C(b)_{k}^{\beta}}) 
=
m \left(
c_{k}^{\beta}
-
\sum_{1 \leq j \leq J, 1 \leq \alpha \leq \alpha(j)} 
q_{jk}^{\alpha \beta} a_{k}^{\beta}
\right)
= 0
\] 
for any $j$ and $\beta$. 
Thus also in this case, $L'\vert_{X_b}$  belongs to $\Pic^0_{X_b}$. 
Therefore $L' \vert_{B'}$ belongs to $\Pic^0_{f^{-1}(B')/B'}$.

We have just shown the lemma when $\Delta$ is irreducible. 
Now we consider a general case, i.e., 
$\Delta = \bigcup_{i=1}^{I_1} \Delta_i$. 
For each $\Delta_i \ (1 \leq i \leq I_1)$, 
take a divisor $Z_i$ of $\Delta_i$ and 
$\sum_{1 \leq i leq I_1, 1 \leq j \leq J_i}e_{ij}\Gamma_{ij}
$ in the same way as above. 
If we set 
\[
B' = B - \left( \vert Z_1 \vert \cup \cdots 
\cup \vert Z_{I_1} \vert \cup (\bigcup_{i,j} 
\vert \Delta_i \vert \cap \vert \Delta_j \vert ) \right),
\] 
then $B'$ is a big open set, and it is easy to see that 
$L^{\otimes m} \otimes \OO_X (- \sum_{ij} e_{ij} \Gamma_{ij}) \vert_{B'}$ 
belongs to $\Pic^0_{f^{-1}(B')/B'}$.
\QED

{\it Proof of Theorem~\ref{theorem:main}}\quad
First we prove the first assertion of the theorem.  
In virtue of Lemma~\ref{lemma:for:main:reduction}, 
by taking a suitable generically finite cover of $B$, 
we may assume that $f : X \to B$ has a section. 
Moreover, by \cite[Theorem 8.2]{Jong}, 
there is a surjective generically finite morphism 
$\widetilde{B} \to B$ of arithmetic varieties 
such that $\widetilde{B}$ is regular. 
Thus, by Lemma~\ref{lemma:for:main:reduction}, 
we may also assume that $B$ is regular. 

We follow the notation of lemma~\ref{lemma:for:main:finite}, 
and let $L^{\otimes m} \otimes \OO_X (- \sum_{ij} e_{ij} \Gamma_{ij})$ 
be a line bundle on $B$ whose restriction to a big open set 
$B'$ of $B$ belongs to $\Pic^0_{f^{-1}(B')/B'}$. 
For simplicity, we set $E = - \sum_{ij} e_{ij} \Gamma_{ij}$. 
Then 
\begin{align*}
\quad & \adeg\left( \acherncl_1(\overline{L^{\otimes m}})^2 \cdot 
\acherncl_1\left(f^*(\overline{H})^{d}\right)\right) \\
& =
\adeg\left( 
\left(\acherncl_1(\overline{L^{\otimes m} \otimes \OO_X(E)}) 
      -\acherncl_1(\overline{\OO_X(E)})\right)^2 
\cdot 
\acherncl_1\left(f^*(\overline{H})\right)^{d}
     \right) \\
& = 
\adeg\left( \acherncl_1(\overline{L^{\otimes m} \otimes \OO_X(E)})^2 
\cdot \acherncl_1\left(f^*(\overline{H})\right)^{d}\right) \\
& \quad -
2 \adeg\left( \acherncl_1(\overline{L^{\otimes m} \otimes \OO_X(E)}) 
\cdot \acherncl_1(\overline{\OO_X(E)}) \cdot
\acherncl_1\left(f^*(\overline{H})\right)^{d}\right) 
+
\adeg\left( \acherncl_1(\overline{\OO_X(E)})^2 \cdot
 \acherncl_1\left(f^*(\overline{H})\right)^{d}\right).
\end{align*}
Since $\deg (L^{\otimes m} \otimes \OO_X(E) \vert_C) = 0$ 
for any vertical curve $C$ lying over $B'$, 
the second term in the last  expression becomes zero. 
Moreover, for the third term in the last expression, we have 
\[
\adeg\left( \acherncl_1(\overline{\OO_X(E)})^2 \cdot
 \acherncl_1\left(f^*(\overline{H})^{d}\right)\right)
=
\sum_{i = 1}^{I_1} \deg_{H}(\Delta_i) \cdot
\left( \sum_{1 \leq j,k \leq J_i} e_{ij} e_{ik} q^i_{jk} \right),
\]
where $q^i_{jk} = \dim_{k(\Delta_i)}(\Gamma_{j,k(\Delta_i)} 
\cap \Gamma_{k,k(\Delta_i)})$. 
From the proof of lemma~\ref{lemma:for:main:finite}, 
this value is non-positive. 
Moreover the equality holds 
if and only if $e_{i1} = \cdots  = e_{iJ_i}$ for $1 \leq i \leq I_1$. 
To sum up, we get
\begin{equation*}
\adeg\left( \acherncl_1(\overline{L^{\otimes m}})^2 \cdot 
\acherncl_1\left(f^*(\overline{H})^{d}\right)\right) 
\leq 
\adeg\left( 
\left(\acherncl_1(\overline{L^{\otimes m} \otimes \OO_X(E)}) 
\right)^2 \cdot \acherncl_1\left(f^*(\overline{H})^{d}\right)
     \right). 
\end{equation*}

Next let $h_L'$ be an admissible line bundle on $L$. 
Then by Lemma~\ref{lemma:for:main:infinite} and Proposition~\ref{prop:main}, 
\begin{multline*}
\adeg\left( 
\left(\acherncl_1(\overline{L^{\otimes m} \otimes \OO_X(E)}) 
\right)^2 \cdot \acherncl_1\left(f^*(\overline{H})\right)^{d}
     \right). \\
\leq 
\adeg\left( 
\left(\acherncl_1(L^{\otimes m} \otimes \OO_X(E), h_L'{}^m) 
\right)^2 \cdot \acherncl_1\left(f^*(\overline{H})\right)^{d}
     \right). \\
=
-2 m^ 2 h_{\OO_{J_K}(\Theta_K)}^{\overline{B}}([L_K]). 
\end{multline*}
Thus we get the first assertion of Theorem~\ref{theorem:main}.
 
Now assuming that $H$ is ample and $\cherncl_1(\overline{H})$ is positive, 
we consider when the equality holds. 

Let $g : \widetilde{B} \to B$ be a surjective generically finite morphism of 
arithmetic varieties such that $\widetilde{B}$ is regular 
and $\widetilde{f} : \widetilde{X} \to \widetilde{B}$ has a section, 
where $\widetilde{X} = X \times_B \widetilde{B}$ and 
\[
\begin{CD}
\widetilde{X} @>{\widetilde{g}}>> X \\
@V{\widetilde{f}}VV @V{f}VV \\
\widetilde{B} @>{g}>> B
\end{CD}
\]
is the induced morphism. 
Let us set $\widetilde{L} = \widetilde{g}^*(L)$ and 
$\widetilde{H} = g^*(H)$. 

Now let us assume the condition (a) and (b) in the second assertion of 
the theorem. By Lemma~\ref{lemma:for:main:finite}
there are a big open set $\widetilde{B'}$ of $\widetilde{B}$, 
a positive integer $m$ and a vertical divisor $\Gamma$ of 
$\widetilde{X}$ such that 
$g \circ \widetilde{f} (\Gamma) \subset B \setminus B''$ 
and that 
$\widetilde{L}^{\otimes m} 
\otimes \OO_{\widetilde{X}}(\Gamma) \left\vert_{\widetilde{B'}} \right.$ 
belongs to 
$\Pic^0_{\widetilde{X}/\widetilde{B}} \left\vert_{\widetilde{B'}} \right.$. 

\begin{Claim}
If $\widetilde{K}$ denotes the function field of $\widetilde{B}$, 
then 
\[
\adeg\left( \acherncl_1(\overline{\widetilde{L}})^2 \cdot 
     \acherncl_1\left(\widetilde{f}^*(\overline{\widetilde{H}})\right)^{d} 
     \right)
=
-2 [\widetilde{K} : K] \  
\widehat{h}_{\OO_{J_{\overline{K}}}
(\Theta_{\overline{K}})}^{\overline{B}}([L_K]). 
\]
\end{Claim}

\Proof
By Proposition~\ref{prop:main}, 
we get
\begin{align*}
\adeg\left( \acherncl_1\left(
     \overline{\widetilde{L}}^{\otimes m} \otimes \OO_{\widetilde{X}}(\Gamma)
                       \right)^2 \cdot 
     \acherncl_1\left(\widetilde{f}^*(\overline{\widetilde{H}})\right)^{d} 
     \right)
& =
-2 m^2 \widehat{h}_{\OO_{J_{\overline{K}}}
(\Theta_{\overline{K}})}^{(\widetilde{B}, g^*(\overline{H}))}([L_K]) \\
& =
-2 m^2 [\widetilde{K} : K] \ \widehat{h}_{\OO_{J_{\overline{K}}}
(\Theta_{\overline{K}})}^{\overline{B}}([L_K])
\end{align*}
On the other hand, since 
\[
\acherncl_1\left(
     \overline{\widetilde{L}}^{\otimes m} \otimes \OO_{\widetilde{X}}(\Gamma)
                       \right)^2
= 
m^2 \acherncl_1(\widetilde{L})^2 
+ 2 m \acherncl_1(\widetilde{L}) 
      \cdot \acherncl_1(\OO_{\widetilde{X}}(\Gamma)) 
+ \acherncl_1(\OO_{\widetilde{X}}(\Gamma))^2
\]
and $\widetilde{f}^*(\widetilde{H}) = f^*(g^*(H))$, 
we get
\[
\adeg\left( \acherncl_1\left(
     \overline{\widetilde{L}}^{\otimes m} \otimes \OO_{\widetilde{X}}(\Gamma)
                       \right)^2 \cdot 
     \acherncl_1\left(\widetilde{f}^*(\overline{\widetilde{H}})\right)^{d} 
     \right)
=
m^2
\adeg\left( \acherncl_1\left(
     \overline{\widetilde{L}}
                       \right)^2 \cdot 
     \acherncl_1\left(\widetilde{f}^*(\overline{\widetilde{H}})\right)^{d} 
     \right).
\]
by projection formula 
(Note that $g \circ \widetilde{f} (\Gamma) \subset B \setminus B''$ ). 
Thus we obtain the claim. 
\QED

From the claim, we get 
\[
\adeg\left(\acherncl_1(\overline{L})^2 \cdot 
           \acherncl_1(f^*(\overline{H}))^d 
     \right)
= 
-2 \widehat{h}_{\OO_{J_{\overline{K}}}
(\Theta_{\overline{K}})}^{\overline{B}}([L_K])
\]
by projection formula. 

Next we assume that 
\[
\adeg\left(\acherncl_1(\overline{L})^2 \cdot 
           \acherncl_1(f^*(\overline{H}))^d 
     \right)
= 
-2 \widehat{h}_{\OO_{J_{\overline{K}}}
(\Theta_{\overline{K}})}^{\overline{B}}([L_K])
\]
Then by projection formula, we have
\[
\adeg\left( \acherncl_1(\overline{\widetilde{L}})^2 \cdot 
     \acherncl_1\left(\widetilde{f}^*(\overline{\widetilde{H}})\right)^{d} 
     \right)
=
-2 \widehat{h}_{\OO_{J_{\overline{K}}}
(\Theta_{\overline{K}})}^{(\widetilde{B}, g^*(\overline{H}))}([L_K]). 
\] 
Let $\widetilde{\Delta}$ be the set of critical values 
of $\widetilde{f}$ 
and $\widetilde{\Delta} = \bigcup_{i=1}^{I} \widetilde{\Delta}_i$ 
be the irreducible decomposition 
of $\widetilde{\Delta}$, where 
$\widetilde{\Delta}_1, \ldots, \widetilde{\Delta}_{I_1}$ 
are divisors on $\widetilde{B}$ such that 
$g(\widetilde{\Delta}_i)$ are also divisors on $B$ for 
$1 \leq i \leq I_1$, 
$\widetilde{\Delta}_{I_1 +1}, \ldots, \widetilde{\Delta}_{I_2}$ 
are divisors on $\widetilde{B}$ such that 
$\codim_B (g(\widetilde{\Delta}_i)) \geq 2$ 
for $I_1 +1 \leq i \leq I_2$, 
and 
$\widetilde{\Delta}_i \ (i \geq I_2)$ satisfy 
$\codim_{\widetilde{B}} (\widetilde{\Delta}_i) \geq 2$. 
Then we take $\sum_{1 \leq i \leq I_2, 1 \leq j \leq J_i} e_{ij} 
\Gamma_{ij}$ as in Lemma~\ref{lemma:for:main:finite} 
(which is applied to $\widetilde{f} : \widetilde{X} \to \widetilde{B}$). 
If we look back closely the proof of the first assertion of 
the theorem, we find that the equality holds 
if and only if (a) $e_{i1} = \cdots = e_{iJ_i}$ for $1 \leq i \leq I_1$ 
and (b) $\overline{L}$ is flat along fibers 
(Note that the reason we need to consider $I_1$ and $I_2$ is that 
$\deg_{g^*(H)}(\widetilde{\Delta}_i) = 0$ for 
$I_1 +1 \leq i \leq I_2$). 
Moreover the condition (a) is equivalent to the existence 
of a big open set $B''$ of $B$ such that $\deg(L\vert_C) = 0$ 
for any fibral curves $C$ lying over $B''$.
This proves the second assertion. 
\QED

\end{document}